\documentclass[11pt]{amsart}

\usepackage{amsmath}
\usepackage{amssymb}
\usepackage{graphicx}

\newtheorem{theorem}{Theorem}[section]

\newtheorem{lemma}[theorem]{Lemma}

\theoremstyle{definition}
\newtheorem{definition}[theorem]{Definition}
\newtheorem{example}[theorem]{Example}
\newtheorem{remark}[theorem]{Remark}

\numberwithin{equation}{section}

\def\bfc{{\bf c}}

\def\bfe{{\bf e}}
\def\bff{{\bf f}}

\def\bfu{{\bf u}}
\def\bfv{{\bf v}}
\def\bfw{{\bf w}}

\def\bfx{{\bf x}}
\def\bfy{{\bf y}}

\def\tri{{\triangle}}

\title[On Locality of Harmonic GBCs]
{On Locality of Harmonic Generalized Barycentric Coordinates 
and Their Application to Solution of the Poisson Equation}
\author[C. Deng]{Chongyang Deng}
\address{School of Science, Hangzhou Dianzi  University, 
Xiasha, Hangzhou 310018, China. 
This author is partially supported by National Nature Science Foundation of China (NSFC) 
under the project numbers: 61872121 and 61761136010.}
\author[M. -J. Lai]{Ming-Jun Lai}
\email{mjlai@uga.edu}
\address{Department of Mathematics,
University of Georgia, Athens, GA 30602.}

\begin{document}

\begin{abstract}
We first extend the construction of generalized barycentric coordinates (GBC) based on the 
vertices on the boundary of a polygon $\Omega$ to a new kind of GBCs based on 
vertices inside the $\Omega$ of interest. 
For clarity, the standard GBCs are called boundary GBCs while the new GBCs are called interior GBCs.  
Then we present an analysis on these two kinds of harmonic GBCs to show that each GBC function 
whose value is $1$ at a vertex (boundary or interior vertex of $\Omega$) decays to zero away from its 
supporting vertex  exponentially fast except for a trivial example. 
Based on the exponential decay property,  we explain how to approximate the harmonic
GBC functions locally. That is, due to the locality of these two kinds of GBCs, 
one can approximate each of
 these GBC functions by its local versions which is supported over a sub-domain of $\Omega$.  
The local version of these GBC function
will help reduce the computational time for shape deformation in graphical design.  
Next, with these two kinds of GBC functions at hand, we can use them to approximate the solution 
of the Dirichlet problem of the Poisson equation.  This may provide a more efficient way
to solve the Poisson equation by using a computer which has graphical processing unit(GPU) 
with thousands or more processes  than the standard methods using a computer with one or few 
CPU kernels.  
\end{abstract}


\maketitle

\section{Introduction}
Generalized barycentric coordinates (GBC) provide a simple and convenient way to
represent a surface over the interior  of a polygon by weighted combinations of 
the control points. They are widely used in computer graphics  and 
related areas. See, e.g. \cite{F03}, \cite{HF06}, \cite{FHK06}, \cite{JSW05}, \cite{KB15}, \cite{HS17}, 
\cite{JBPS11}, \cite{Zetc14}, \cite{CXC2018}, and etc. 
In addition, they found their applications in numerical solution of 
partial differential equations.
We refer the interested reader to \cite{MRS14}, \cite{RGB14}, \cite{FL16}, \cite{LL17} and etc. 
The study of generalized barycentric coordinates 
(GBC) started from a seminal work in \cite{W75}. Since then, there are many
GBCs which have been constructed. See  \cite{F15}  and \cite{A18} for 
numerous GBC studied in literature.

However, many of them are defined as a combination of all
control points. So changing a single control point will lead to a change to the surface entirely 
inside the polygon. This is not convenient for practical applications. 
For example, in shape and image deformation, a control point should 
only influence the surface nearby the places where the changes are made. 
This is the so-called locality property which is extremely important to a designer.
In additional, without locality, the geometry design requires a large memory consumption of 
all GBC functions as a polygon with $n$ sides may have a very large $n$ for a reasonable object. 
Such a design will be computationally expensive. 
To overcome this difficulty of computation, several approaches have been developed. 
For example, Zhang et al. \cite{Zetc14} proposed a family of local barycentric coordinates (LBC) 
through a convex optimization approach. 
For another example, subdividing barycentric coordinates (SBC) \cite{ADH16} 
inherit local property from local support of subdivision surfaces. Another example is 
blended barycentric coordinates (BBC) \cite{APH17} which are blended from mean value 
coordinates over 
the triangles of the constrained Delaunay triangulation 
of the input polygon which possesses a locality. In fact, 
there is a folklore in the community of computer graphics that harmonic GBC functions have 
desired intuitive locality behavior in all situations during a geometry design 
although the maximum principle shows that any harmonic GBC function is nonnegative even over 
a polygon which is strongly 
concave and never be zero inside the polygon.  However, even though the researchers and 
practitioners in the community have observed this locality for a long time,  
there is no mathematical justification of this locality in the literature. 
Nevertheless there is a trivial example that 
harmonic GBCs over a rectangular domain do not have any local property. 
This makes a doubt that such a locality can be established for polygons with more than 4 sides.   

In this paper, our first goal is to explain the 
nice locality of harmonic GBC function over a polygon which is more complicated 
than a rectangle mathematically. To this end, let us define an exponential decay property.  

\begin{definition}
Let $f(\bfx)$ be a function defined on a polygon $\Omega\subset \mathbb{R}^d, d\ge 2$ 
and suppose that $f(\bfx_0)=1$  for $\bfx_0\in \Omega$. 
$f$ has an exponentially decay property away from its supporting vertex 
$\bfx_0$ if $|f(\bfx)|\le C \exp(- c\|\bfx-\bfx_0\|)$ for positive constants $C$ and $c$. 
Such a property is called e-locality for short.  
\end{definition}

That is,  we shall  discuss the e-locality of harmonic GBC functions in this paper.  
Mainly, we present a method  to analyze that each  harmonic coordinate with $1$ at a vertex of 
the polygon $\Omega$ decays away from its vertex exponentially over $\Omega$. We have
to exclude the case when $\Omega$ is a parallelogram as the harmonic GBCs are 
simply bilinear functions.  We will 
discuss this pathological example more later in this paper. This makes sense as for most 
applications,$\Omega$ of interest usually has a lot of  sides.  
Thus, any change of the control point  of a harmonic coordinate
function (the value at the supporting vertex) will affect the surface nearby and the change will 
decay exponentially to 0 away from the control point over the polygon. 
Next let us be more precise about the exponential decay property as follows.
 
To do so, let us introduce harmonic GBC functions now. 
Given a polyhedron $\Omega=P_N\in \mathbb{R}^d, d\ge 2$ of $N$  
vertices $\bfv_i, i=1, \cdots, N$, we say
$P_N$ is an admissible polyhedron if $P_N$ admits a simplicial partition $\triangle$ with 
these vertices of $P_N$ being 
the vertices of $\triangle$. In $\mathbb{R}^2$, any polygon $P_N$ is admissible. Note that in 
$\mathbb{R}^3$, one can
have a polyhedron which does not admit a simplicial partition without adding more vertices.  
The admissible condition allows us 
to define a piecewise linear function $\ell_i$  associated with each vertex $\bfv_i$ on the 
boundary of $\Omega$  satisfying  $\ell_i(\bfx)=1$ if $\bfx=\bfv_i$ and $0$ 
if $\bfx=\bfv_j, j\not=i$  and $\ell_i(\bfx)$ is linear on 
each lower dimensional boundary simplex of $\Omega$.  For each $\ell_i$, let $\phi_i$ 
be the function solving the following minimization problem:
\begin{equation}
\label{vPDE}
\begin{cases}
\min_{u\in C^\infty(\Omega)} &\int_{\Omega}|\nabla u|^2, \cr
\qquad u &= \ell_i, \quad \bfx\in \partial \Omega,
\end{cases}
\end{equation}
where $\nabla$ is the gradient operator. It is known that the minimizer $\phi_i$ satisfies 
the following Laplace equation:
\begin{equation}
\label{PDE}
\begin{cases}
\Delta u &=0, \quad \bfx\in \Omega\cr
u &= \ell_i, \quad \bfx\in \partial \Omega,
\end{cases}
\end{equation}
where $\Delta$ is the standard Laplace operator. 
That is,  the minimizer $\phi_i$ is a harmonic function.   It is easy to see that these $\phi_i$ satisfy the three 
properties in (\ref{GBC}). In fact, one can use (\ref{PDE}) to verify these properties. 
\begin{equation}
\label{GBC}
\begin{cases}
\displaystyle \sum_{i=1}^n \phi_i(\bfx)  &=  1 \quad \bfx\in P_N\cr
\displaystyle \sum_{i=1}^n \phi_i(\bfx)\bfv_i &=  \bfx \quad \bfx\in P_N\cr
\phi_i(\bfx) &\ge 0, \quad i=1, \cdots, N.
\end{cases}
\end{equation}
Hence, these $\{\phi_1, \cdots, \phi_N \}$ are  generalized barycentric coordinates (GBC) and  
these are called harmonic GBCs (cf. \cite{F15}).


These GBCs are based on the corners of $\Omega$. We now extend the GBC functions 
to more boundary points of 
$\Omega$.  Let us add more points on the boundary and inside of $\Omega$ and 
let $\triangle$ be a triangulation of $\Omega$ if $\Omega$ is a polygonal domain in $\mathbb{R}^2$. For $d> 2$, 
$\triangle$ is a simplicial partition of $\Omega$ when $\Omega$ is an admissible polygon 
in $\mathbb{R}^d$.  
Let $V_B=\{\bfv\in \triangle: \bfv\in 
\partial \Omega\}$ be the set of boundary vertices of $\Omega$ and 
$V_I= \{\bfv\in \triangle: \bfv\in \Omega^\circ\}$
be the collection of the interior vertices of $\Omega$, where $\triangle^\circ$ stands for the 
interior of $\triangle$.  
As $V_B$ contains more than the corners of $\Omega$,   
we can view $\Omega$ is a degenerated admissible polyhedron with vertices in $V_B$.  
Let $\phi_1, \cdots, \phi_N$ be  GBC functions associated with the vertices of $V_B$ 
defined before.  

The second goal of this paper is to introduce another kind 
of GBC functions which are based on interior points of $\Omega$ as follows.  
For $\triangle$ of $\Omega$, 
let $S^0_1(\triangle)$ be the continuous piecewise linear functions over 
simplicial partition $\triangle$.  That is, 
$S^0_1(\triangle)$ is the space of continuous linear splines.  
For each $\bfv_i\in V_I= \{\bfv\in \triangle: \bfv\in \Omega^\circ\}$, 
let $h_i\in S^0_1(\triangle)$ be the hat function
which is a continuous piecewise linear function satisfying $h_i(\bfv_j)=\delta_{ij}$, for all 
$\bfv_j\in V_I$. For each $h_i$, let $\psi_i\in H^\infty(\Omega)$ be the solution to the following
boundary value problem:
\begin{equation}
\label{PDE2}
\begin{cases}
\Delta \psi_i &= h_i, \quad \bfx\in \Omega\cr
\psi_i &= 0, \quad \bfx\in \partial \Omega,
\end{cases}
\end{equation}
for $i=1, \cdots, M$, where $M=\#(V_I)$. 
It is easy to see that  there exist such functions $\psi_i$ by solving the Poisson 
equation with zero boundary condition. Let us further explain their properties. 
A triangle $T\in \triangle$ is a boundary triangle if $T$ has one
side on the boundary of $\Omega$ or a vertex at the boundary of $\Omega$.
Let $\triangle^\circ$ be the union of all simplexes in $\triangle$ without boundary triangles. 
 
\begin{figure}[htbp]
\centering
\includegraphics[width = 0.8\textwidth]{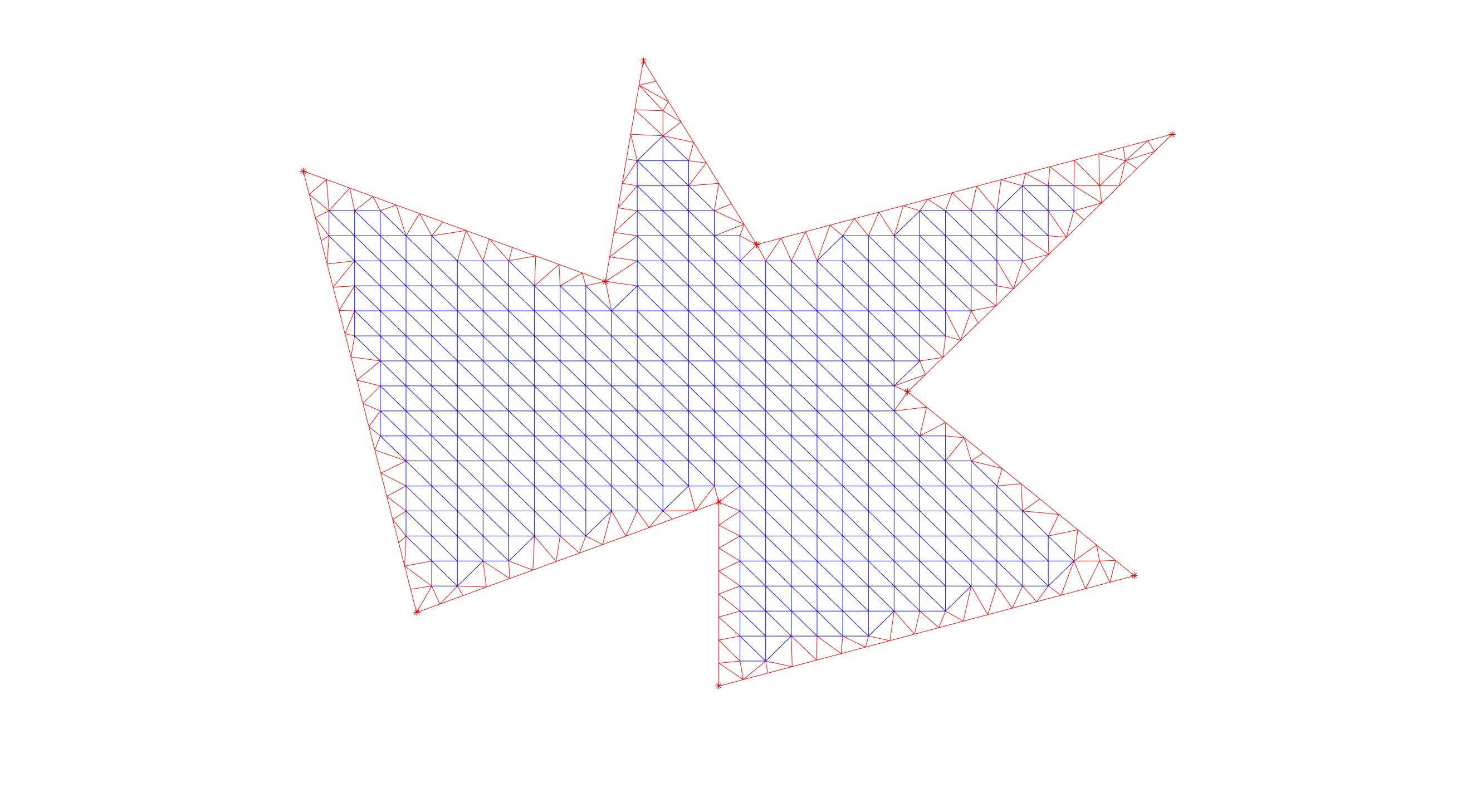}
\caption{Illustration of $\triangle^\circ$ (in blue) inside a triangulation $\triangle$ 
(all triangles in red and blue) \label{tri0}}
\end{figure}

It is easy to see that $\psi_i$'s  satisfy the properties in (\ref{GBC2}) below 
and hence, $\psi_i$ are called interior GBC functions.  
\begin{lemma}
\label{lem1}
Let $\psi_i, i=1, \cdots, M$ be functions defined above. Then 
\begin{equation}
\label{GBC2}
\begin{cases}
\displaystyle 
\sum_{i=1}^m \Delta \psi_i(\bfx)  &=  1, \quad \bfx\in \triangle^\circ \cr
\displaystyle 
\sum_{i=1}^m \Delta \psi_i(\bfx)\bfv_i &=  \bfx, \quad \bfx \in \triangle^\circ\cr
\psi_i(\bfx) &= 0, \quad \bfx\in \partial \Omega,
\end{cases}
\end{equation}
where $\bfv_i, i=1, \cdots, M$ are interior vertices of $\triangle$.  See Figure~\ref{tri0} for an illustration of
interior domain $\triangle^\circ$. 
\end{lemma}
\begin{proof}
Since $\sum_{i=1}^m h_i(\bfx)=1$  and $\sum_{i=1}^m \bfv_i h_i(\bfx)= \bfx$ over $\triangle^\circ$, we see that
$\sum_{i=1}^m \Delta \psi_i(\bfx) =  \sum_{i=1}^m h_i(\bfx)=1$ for $\bfx\in \triangle^\circ$. Similarly, we have
the second equation in (\ref{GBC2}). The third equation follows from the boundary condition of $\psi_i$.  
\end{proof}

It is known that there is no analytic representation of harmonic GBCs $\phi_i$ when 
$\Omega$ is not a parallelogram.  One has to compute them numerically.  
A standard approach is to use finite element
method or more generally, multivariate spline method (cf. \cite{DFL18}). Similarly, $\psi_i$
have no analytic representation to the best of the authors' knowledge.  Of course, one can solve
the Poisson equation in (\ref{PDE2}) numerically to obtain an approximation of $\psi_i$.  Indeed, 
let $\triangle_k$ be the $k$th uniform refinement of $\triangle$ for an integer $k\ge 1$ 
if $d=2$ or $\triangle_k$ be a refined 
simplicial partition of $\triangle$ if $d\ge 3$ with $|\triangle_k|< |\triangle|$, where 
$|\triangle|$  the largest diameter of simplexes in $\triangle$. Thus, 
$|\triangle|$ is called the size of $\triangle$.    
Let $S_i=S_{\phi_i}$ be a spline approximation of
$\phi_i$ over triangulation $\triangle_k$ and 
similarly $R_i=S_{\psi_i}$ be a spline approximation of $\psi_i$ over $\triangle_k$.  

One of the significances of these
functions is that we can use $\phi_i, i=1, \cdots, N$ and $\psi_j, j=1,\cdots, M$ to approximate 
the solution of any Dirichlet boundary value problem of the Poisson equation:
\begin{equation}
\label{Poisson}
\begin{cases}
-\Delta u &=f, \quad \bfx\in \Omega\cr 
\quad u &= g, \quad \bfx\in \partial \Omega,
\end{cases}
\end{equation}
for any given $f\in L^2(\Omega)$ and $g\in L^2(\partial \Omega)$. 
For simplicial partition  $\triangle$ of $\Omega$, 
we are interested in the approximation of $S_i$ and $R_i$ when $|\triangle|$ is small.  
Indeed, the following is one of the main results in this paper: 
\begin{theorem}
\label{thm1}
Suppose that $f$ is a piecewise continuous function in $L^2(\Omega)$ and $g$ is a piecewise
continuous function in $L^2(\partial\Omega)$.  
Let $\triangle$ be a simplicial partition of $\Omega$
such that  $H_f(x,y):=\sum_{i=1}^m f(\bfv_i) h_i(x,y)$
is a linear spline approximation of $f$ over $\Omega$, e.g. 
$f(x,y)- H_f(x,y)= O(|\triangle|^2)$ and 
$L_g:= \sum_{\bfv_i\in \partial\Omega} g(\bfv_i)S_i$ is a linear spline 
approximation of $g$ on $\partial\Omega$, i.e., $\|g- L_g\|_{L^2(\partial \Omega)} 
= O(|\triangle|^2)$. 
Then the solution $u$ to (\ref{Poisson}) can be approximated by 
\begin{equation}
\label{DL19}
u \approx \sum_{\bfv_i\in \triangle^\circ} f(\bfv_i) \psi_i + 
\sum_{\bfv_i\in \partial\Omega} g(\bfv_i)\phi_i \approx  
\sum_{\bfv_i\in \triangle^\circ} f(\bfv_i) R_i + 
\sum_{\bfv_i\in \partial\Omega} g(\bfv_i)S_i. 
\end{equation}
That is, letting $L_u= L_f+ L_g$ with $L_f= \sum_{\bfv_i\in \triangle^\circ} f(\bfv_i) R_i$ and
$L_g=\sum_{\bfv_i\in \partial\Omega} g(\bfv_i)S_i$, we have 
$$\|\nabla( u-L_u)\|=  O(|\triangle|).$$ 
\end{theorem}
\begin{proof}
As expected, $L_u$ is just like the standard finite element solution.
We shall give a proof in a later section. See \S 3 for a detailed proof and numerical 
experimental results which 
support the statement in (\ref{DL19}). In addition, a comparison  with the numerical results 
from the standard finite element method will be shown. 
\end{proof}

\begin{remark}
We remark that the GBC solution $L_u= L_f+ L_g$ is not an FEM solution. For simplicity, let us say the solution 
$u$ satisfying the zero boundary condition. In this case, $L_u= L_f$ which is a linear combination of $f(\bfv_i), 
\bfv_i\in \tri^\circ$ while the FEM solution is a linear combination of the coefficient vector $\bfc$ 
which is the solution
of the linear system $K\bfc= M\bff$, where $K$ and $M$ are the stiffness and mass matrices, respectively.  
\end{remark}

Let us continue to discuss the locality of these harmonic GBC functions. First of all, 
we need to explain 
that a numerical harmonic GBC function $S_i$  approximates the exact harmonic GBC function 
$\phi_i$ very well. Similar for $R_i$ for $\psi_i$.   
Therefore, the e-locality of $\phi_i$ and $\psi_i$ can be seen from the e-locality of $S_i$
and $R_i$, respectively.  To avoid a pathological example,
we assume that the number of distinct boundary vertices is more than $4$ as most applied 
problems have a lot of boundary vertices.  

For a simplicial partition $\triangle$  of $\Omega$, we say it is $\beta$-quasi-uniform if 
there is a positive number $\beta>0$ such that 
\begin{equation}
\label{beta} 
\max_{T\in \triangle} \frac{|\tri|}{\rho_T} \le \beta, 
\end{equation}
where  $\rho_T$ is the radius of the inscribed ball of simplex $T$.  
Let $S^r_n(\triangle)$ be the spline space of degree $d$
and smoothness $r\ge 0$, e.g. when $n=1$ and $r=0$, $S^0_1(\triangle)$ is the 
standard continuous finite element space. 
In general, we can use any $n\ge 1$ and $r\ge 0$ as long as $n\ge 3r+2$ if 
$\Omega\subset \mathbb{R}^2$. Otherwise, we
can use $n=1$ and $r=0$ or $n\ge 9$ and $r=1$. See \cite{LS07} for trivariate splines.    
There are two numerical implementation methods  of multivariate splines which can be found in 
\cite{ALW06} and
\cite{S15} as well as other efficient implementations in \cite{M19} and \cite{X19}.  
Recall $\phi_i\in H^1(\Omega)$ be a GBC function in the Sobolev space
$H^1(\Omega)$ satisfying (\ref{vPDE})  
and $S_i$ is the minimizer of the following minimization (\ref{vPDE2}) 
over the spline space $S^r_n(\triangle)$ of degree $n$ and smoothness $r\ge 0$:
\begin{equation}
\label{vPDE2}
\begin{cases}
\min_{u\in S^r_n(\triangle)} &\int_{\Omega}|\nabla u|^2, \cr
 u &= \ell_i, \quad \bfx\in \partial \Omega.
\end{cases}
\end{equation}
We split $S_i$ into $S_{i,0}$ and $G_i$. That is, $S_i= S_{i,0}+G_i\in S^r_n(\triangle)$, 
where $G_i\in S^r_n(\triangle)$ satisfies the boundary condition $G_i=\ell_i$ on 
$\partial \Omega$ and 
$S_{i,0}\in H^1_0(\Omega)\cap S^r_n(\triangle)$.   
Then $S_i$ satisfies the following weak formulation:
\begin{equation}
\label{weak}
\langle \nabla S_{i,0}, \nabla \psi\rangle = -\langle \nabla G_i, \nabla\psi\rangle, 
\quad \forall \psi\in S^r_n(\tri) \cap H^1_0(\Omega). 
\end{equation}
It is easy to see that $S_i$ is a numerical harmonic GBC which approximates $\phi_i$ very well 
in the following sense:
\begin{equation}
\label{cea}
\|\nabla (S_i -\phi_i)\| \le C_{\phi_i} |\triangle|^{d}
\end{equation}
by using the well-known Ce\'a lemma, where $C_{\phi_i}$ is a positive constant dependent on 
$\phi_i$, d, and $\Omega$. 
Furthermore, there is a maximum norm estimate, i.e.
\begin{equation}
\label{maxnormest}
\| S_i - \phi_i\|_{\infty} \le C_{\phi_i} log(|\triangle|) |\triangle|^{d+1}.
\end{equation}
We refer to \cite{C78} for detail. 
Similarly, 
$R_i\in S^r_n(\triangle)$ is the weak solution satisfying
\begin{equation}
\label{weak2}
\langle \nabla R_i, \nabla \psi\rangle = -\langle h_i, \psi\rangle, \quad \forall \psi\in S^r_n(\tri)\cap H^1_0(\Omega). 
\end{equation}
The standard finite element theory shows that $R_i$ satisfies the same inequalities as
$S_i$ in (\ref{cea}) and (\ref{maxnormest}). 
Hence, one can see that the locality of $\psi_i$ can be estimated based on the locality of $R_i$. 
We leave the proof of (\ref{cea}) and (\ref{maxnormest}) to  Appendix. 
 

 
In terms of spline functions, the exponential decay can be recast more precisely. Let
 $\hbox{star}^1(\bfv_i)$ be the union of all simplexes in $\triangle$ sharing vertex $\bfv_i$ and $\hbox{star}^{k}(\bfv_i)$ is the union of all simplexes 
in $\triangle$ sharing vertices in $\hbox{star}^{k-1}(\bfv_i)$ for $k\ge 2$.
We will show that there exists a constant $\sigma\in (0, 1)$ such that 
\begin{equation}
\label{sigma}
|S_i(\bfv)|\le C \sigma^k \hbox{ and } |R_i(\bfv)|\le C \sigma^k,  \hbox{ if } \bfv\not\in 
\hbox{star}^k(\bfv_i)
\end{equation}
for a positive constant $C$ independent of $i$. 
That is, if $\bfv$ is a far away from $\bfv_i$ according to the triangulation $\triangle$ (cf. Definition~\ref{trimeasure}), 
$|S_i(\bfv)|$ is close to zero.  Similar for $R_i$. 
See the statement of Theorem~\ref{main} in the next section for more detail.  

The rest of the paper is devoted to the estimate (\ref{sigma}) 
of the e-locality of $S_i$ and $R_i$ by using a theoretical approach which was used in study of 
the domain decomposition methods for scattered data interpolation 
and fitting (cf. \cite{LS09}) and the convergence 
of discrete least squares (cf. \cite{GS02}).  
We shall present some numerical e-locality of our GBC functions. See \S 4.   
In addition, an application to numerical solution to Poisson equations using our GBC functions 
will be demonstrated in \S 5, where Theorem~\ref{thm1} will be proved. 



\section{The Exponential Decay Property of Boundary and Interior GBC Functions}
\begin{figure}[htbp]
\centering
\includegraphics[width = 0.4\textwidth]{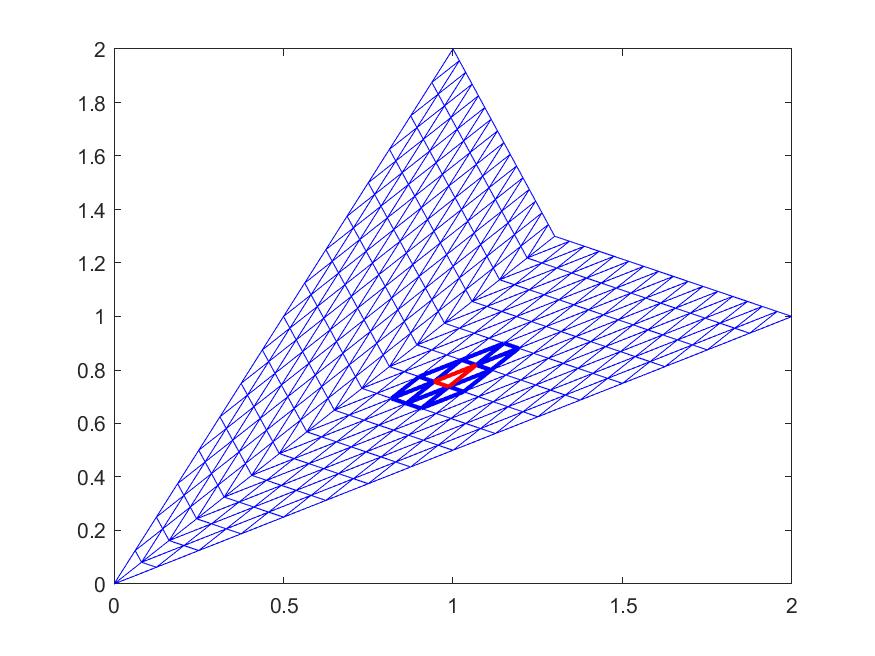}
\includegraphics[width = 0.4\textwidth]{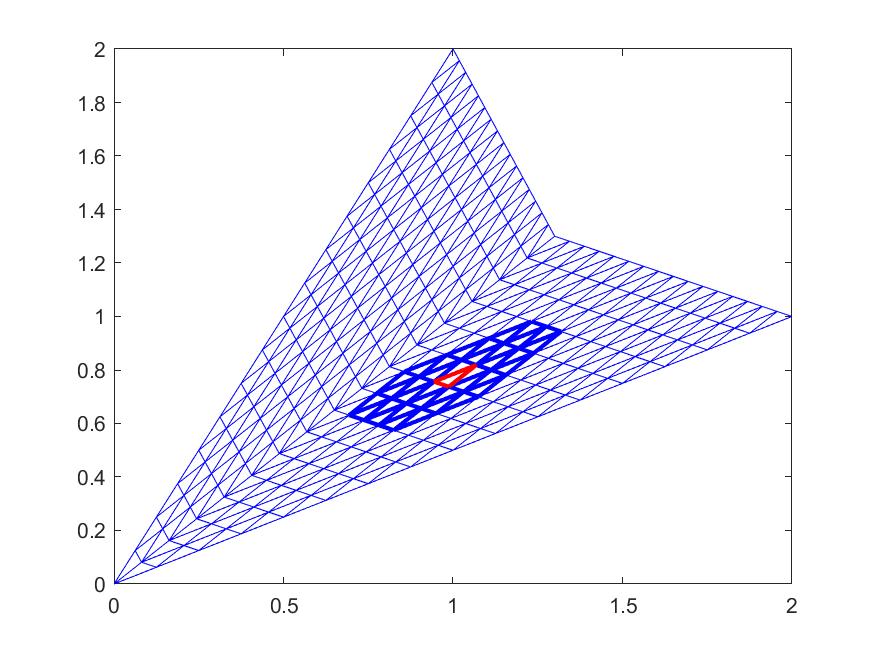}
\includegraphics[width = 0.4\textwidth]{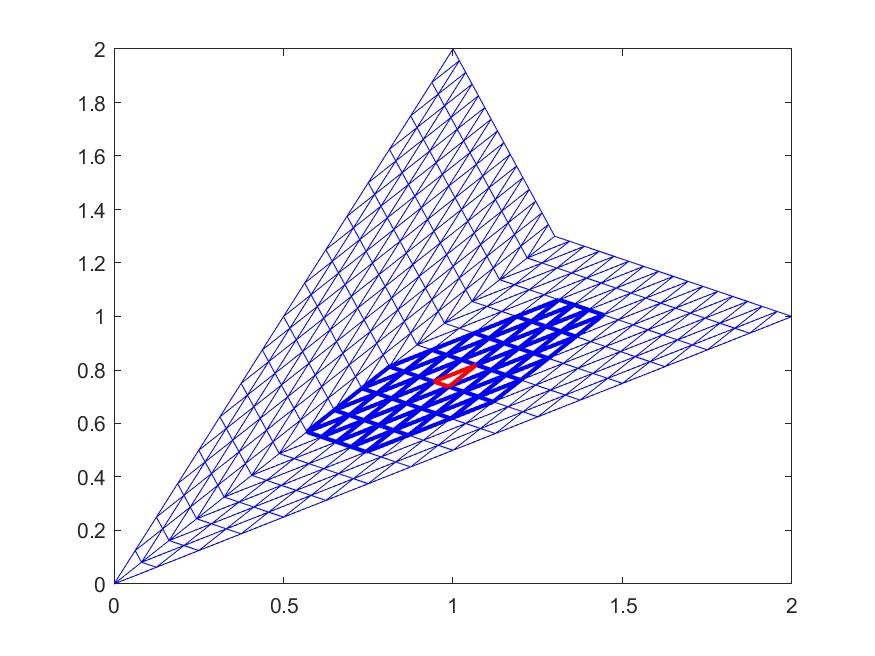}
\caption{Illustration of $\hbox{star}^1(T), \hbox{star}^2(T), \hbox{star}^3(T)$ over
a triangulation with $T$ shown in red  \label{starfig}}
\end{figure}

Let us start with an explanation of some useful concepts,  notations, and definitions on spline 
spaces. For each vertex $\bfv\in \triangle$, 
let $\hbox{star}(\bfv)$ be the collection of all simplexes  from $\triangle$ attached to $\bfv$.
Similarly, for each simplex $T\in \triangle$, let $\hbox{star}(T)$ be the collection of all 
simplexes in $\tri$ connected to $T$. Next for each integer $\ell>1$, let 
$\hbox{star}^\ell(\bfv)$ the the collection of all simplexes in 
$\tri$ which is connected to $\hbox{star}^{\ell-1}(\bfv)$ with $\hbox{star}^1(\bfv)
:= \hbox{star}(\bfv)$ for $\ell>1$. 
Similar for  $\hbox{star}^\ell(T)$. See an example in Figure~\ref{starfig} 
for $\hbox{star}^\ell(T)$ for $\ell=1,2,3$ in the setting of $\mathbb{R}^2$.

Now we define a measure between two points $\bfx, \bfy\in \Omega$ based on triangulation $\triangle$ in the following
sense:
\begin{definition}
\label{trimeasure}
We say a point $\bfx\in \Omega$ is $k$-simplexes away from another point $\bfy\in T \in \triangle$ 
is $\bfx\not\in \hbox{star}^k(T)$, but in $\hbox{star}^{k+1}(T)$.  We say $\bfx$ is far away from $\bfy\in T$ or from $T$ 
is $\bfx\not\in \hbox{star}^k(T)$ for positive integer $k\gg 1$.  
\end{definition}

Next let $B_\xi, \xi\in \mathcal{M}$ be a basis for $S^r_n(\triangle)$, where  $\mathcal{M}$ is an index set. 
For an integer $\ell\ge 1$, we say $B_\xi, \xi\in \mathcal{M}$ 
 is stable and $\ell$-local in the sense that there is a positive constant $K_1$ such that 
for all indices $\xi \in \mathcal{M}$,
\begin{equation}
\|B_\xi\|_{\infty, \Omega} \le K_1, \label{bnd} 
\end{equation}
and there is an integer $\ell\ge 1$  for all $\xi\in \mathcal{M}$, 
\begin{equation} 
\hbox{supp} ( B_\xi ) \subseteq \hbox{star}^\ell(T_\xi), \label{loc} 
\end{equation}
for a simplex $T_\xi\in \tri$ associated with index $\xi$, where the constant $K_1$ may 
be dependent only on $\ell$ and the smallest angle in $\tri$.

There are many spaces which have  stable local bases.  For example, in Euclidean 
space $\mathbb{R}^d$ with $d\ge 2$, 
the continuous spline spaces $S^0_n(\tri)$ have stable local bases with $\ell = 1$ for 
any degree $n\ge 1$.
In $\mathbb{R}^2$  the same is true for the superspline spaces $\mathcal{S}^{r,2r}_{4r+1}(\tri)$
for all $r \ge 1$ of degree $4r+1$ and smoothness $r$ and supersmoothness $2r$ at vertices. 
E.g. $S^{1,2}_5(\triangle)$, the Argyris $C^1$ quintic finite element space with super 
smoothness $C^2$ at vertices.  
In \cite{LS07},  several families of macro-element
spaces defined for all $r \ge 1$ with the same property are explained. 
We refer to \cite{GS02} for more such spaces. 

\begin{figure}[htbp]
\centering
\includegraphics[width = 0.8\textwidth]{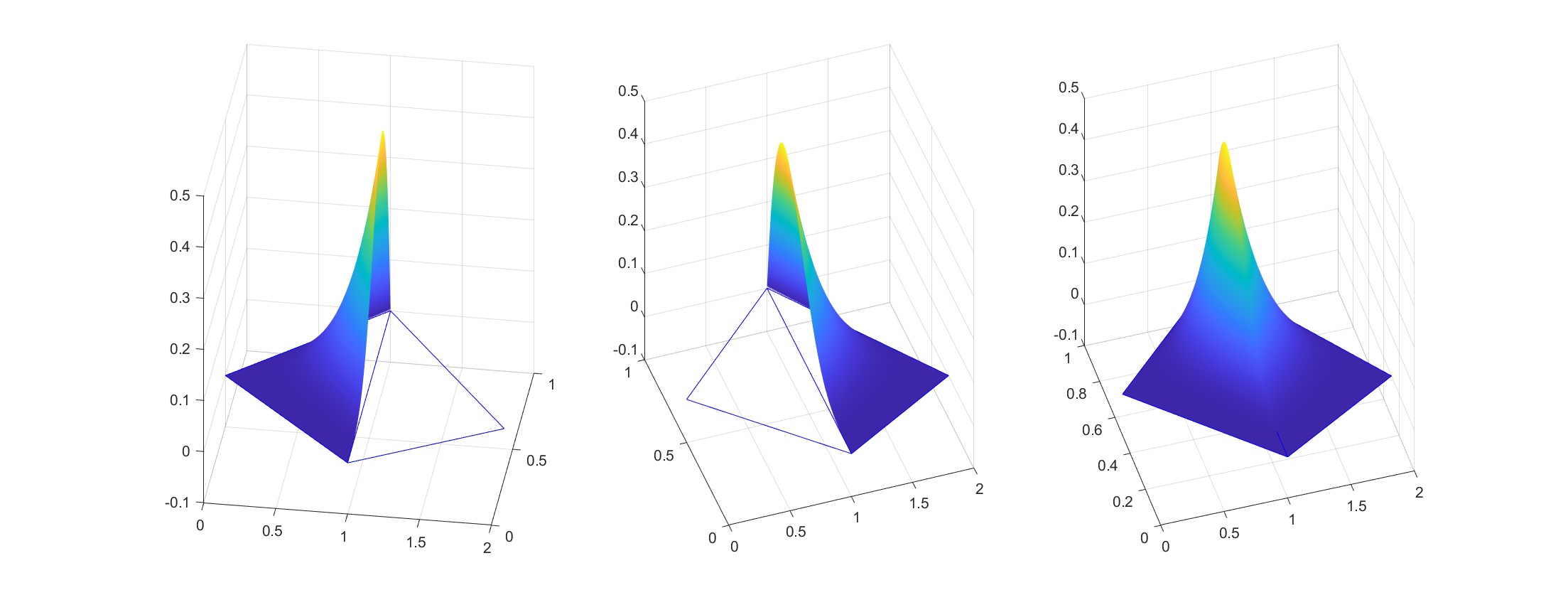}
\caption{Construction of $B_\xi$ over
$T_1\cup T_2$ shown on the right based on $\phi_\xi$ (left) and $\phi_\eta$(middel) \label{Bxi}}
\end{figure}

We now use the format of multivariate Bernstein polynomials to write each polynomial over a 
simplex $T$ (see \cite{CL91}, \cite{B87}, and \cite{LS07} for detail).  For polynomial $u\in \mathbb{P}_n$, we write 
\begin{equation}
\label{Bform}
u(\bfx) =  \sum_{i_0+i_1+\cdots+i_d=n} c_{i_0, \cdots, i_d} \frac{n!}{\prod_{j=0}^d i_j!} \prod_{j=0}^d b_j^{i_j}, 
\end{equation}
where $b_0, \cdots, b_d $ are barycentric coordinates of point $\bfx$ with respect to $T=\langle \bfv_0, 
\cdots, \bfv_d\rangle$. Let 
$\bfc=(c_{i_0,\cdots, i_d}, i_0+\cdots+i_d=n)$ be the coefficient vector of $u$.  Following the proof of Theorem
~2.7 for the $\mathbb{R}^2$ setting and/or the proof of Theorem 15.9 (for the $\mathbb{R}^3$ setting) in \cite{LS07}, we can show 
\begin{equation}
\label{stability2}
\frac{V_T^{1/2}}{K} \|{\bf c}\|_2 \le \|u\|_{2,T} \le V_T^{1/2} \|{\bf c}\|_2 , \quad \forall u\in \mathbb{P}_n 
\end{equation}
for a positive constant $K$ dependent only on $n$, where $V_T$ is the volume of simplex $T$. Letting 
\begin{equation}
\label{Bernstein}
\phi_\xi = \frac{n!}{\prod_{j=0}^d i_j!} \prod_{j=0}^d b_j^{i_j}, 
\end{equation}
with $\xi= \{i_0, \cdots, i_d, T\}$, we note that there are redundancies  
for the collection $\{\phi_\xi, \xi=(i_0, \cdots, i_d, T), T\in \triangle, i_0 
+ \cdots + i_d=n\}$ due to the neighboring simplexes. 
Indeed, if two simplexes $T_1$ and $T_2$ sharing a common face, there is $\phi_\xi$ for 
some $\xi$ associated with $T_1$ and 
$\phi_\eta$ for $\eta$ associated with $T_2$ which are the same function on the common 
face $T_1\cap T_2$. For these two functions, we let $B_\xi$ be the function which is piecewise defined on the 
union $T_1\cup T_2$ with $B_\xi=\phi_\xi$ on $T_1$ and $B_\xi=\phi_\eta$ on $T_1$. Let us delete the index 
$\eta$ from the whole index set  
$\{\xi=(i_0, \cdots, i_d, T), T\in \triangle, i_0 + \cdots + i_d=n\}$.     For simplicity, let us show $B_\xi$ in the
2D setting in Figure~\ref{Bxi}.

For another case, if a vertex $\bfv$ shared by more than one simplex,  there will be 
more than one function $\phi_\xi$ 
with $\xi$ on $T_1, \cdots, \eta$ on $T_k$ which have the common value at the vertex, we 
let $B_\xi$ be the spline function 
which is $\phi_\xi$ on $T_1$, $\cdots, $ $B_\xi=\phi_\eta$ on $T_k$ if $T_1, \cdots, T_k$ 
share the vertex $\bfv$. We delete
all the indices $\eta$ on $T_2, \cdots, T_k$ from the whole index set.   
We do the same thing for other common facets such as edges, ..., $(d-2)$ facets sharing by 
more than one simplex from $\triangle$. 
For the remaining indices  $\xi$, we let $B_\xi=\phi_\xi$.  Let us put these $B_\xi$ together to 
form $\{B_\xi, \xi\in \mathcal{M}\}$ which is a basis for $S^0_n(\triangle)$, where $\mathcal{M}\subset 
\{\xi: T\in \triangle, i_0+\cdots+i_d=n\}$ is the subset  after deleting the redundant 
indices mentioned above.

We are now ready to explain the exponential decay property of the GBC functions. 
Let us first begin with a key result, Lemma~\ref{thm:LS09} whose proof can be found 
in \cite{LS09}. The similar results hold in the multivariate setting which will be 
stated in Theorem~\ref{main1}.   
The proof of Theorem~\ref{main1} relies on the following stability of spline functions.

\begin{lemma}
\label{lai2020}
Suppose that $S^0_n(\triangle)$ is a continuous spline space of degree $n$ 
defined on a $\beta$-quasi-uniform simplicial partition $\triangle$ of polyhedral domain 
$\Omega\subset \mathbb{R}^d$ with $d\ge 2$. Let  $\mathcal{W}=H^1_0(\Omega)\cap S^0_n(\triangle)$ 
be a subspace  with inner product $\langle f, g\rangle_{\mathcal{W}} = \langle \nabla f, 
\nabla g\rangle$ and norm $\|f\|_{\mathcal{W}}:= \|\nabla f\|$. 
Then there exist two positive constants $C_1, C_2$, the following equalities 
\begin{equation}
C_1  \sum_{\xi\in\mathcal{M}} |c_\xi|^2
\le \big \|\sum_{\xi\in \mathcal{M}} c_\xi \mathcal{B}_\xi \big\|_{\mathcal{W}}^2
   \le C_2  \sum_{\xi\in \mathcal{M}} |c_\xi|^2  \label{Wriesz}
\end{equation}
hold for all coefficient vectors $\bfc := \{c_\xi\}_{\xi\in\mathcal{M}}$, 
where  $\{\mathcal{B}_\xi\}_{\xi\in \mathcal{M}}$ is an 1-local basis for $\mathcal{W}$. 
\end{lemma}
\begin{proof}
Let $\mathbb{P}_n$ be the space of polynomials of degree $\le n$ in $\mathbb{R}^d$. 
Let $\mathbb{T}$ be the set of all simplexes $T$ with one vertex at ${\bf 0}$ and $1/\beta 
\le \rho_T \le |T| = 1$, where $\rho_T$ is the radius of the ball inscribed in $T$ and 
$|T|$ is the diameter of the minimum ball containing $T$. Let
$$ C_1 := \inf_{T =\langle \bfv_1, \cdots, \bfv_{d+1}\rangle \in \mathbb{T}} \{ \|\nabla p\|_{T}^2 \hbox{ such that }  
p\in \mathbb{P}_n, \int_T p^2 = 1,   p(\bfv_1)= 0 \}. 
$$
Then there exist sequences $p_k, T_k$ of polynomials and simplexes, 
respectively, such that $p_k \to p\in \mathbb{P}_n$ and $T_k \to T \in \mathbb{T}$ with $\int_T p^2 = 1$ 
and $C_1 = \|\nabla p\|_{T}^2$.   We claim that $C_1 > 0$. Indeed, if
$\|\nabla p\|_{T}^2 = 0$, then $p \in \mathbb{P}_0$.  But then 
using the fact that $p$ vanishes at one vertex, it follows
that $p \equiv 0$, contradicting $\int_T p^2 = 1$.
We have shown that $C_1> 0$ and that it certainly depends only on $\beta$, $n$, and $d$. 

Next let
$$
C_2 := \sup_{T =\langle \bfv_1, \cdots, \bfv_{d+1}\rangle \in \mathbb{T}} 
\{ \|\nabla p\|_{T}^2  \hbox{ such that }  
p \in \mathbb{P}_n, \int_T p^2 = 1, \hbox{ and } 
  p(\bfv_1)= 0 \}. 
$$
Clearly, by using the Markov inequality, it is easy to see that $C_2 < \infty$. 
Certainly,  $C_2$ depends only on $\beta$, $d$, and $n$.   
Now if $T$ is an arbitrary simplex in $\triangle$, then
after translating one vertex to ${\bf 0}\in \mathbb{R}^d$ and
substituting $\bfx = |T|\widetilde{\bfx}$, we see that
for any $u \in \mathcal{W}$,
\begin{equation}
\label{stability} 
C_1 \int_T u^2 \le |T|^{d} \|\nabla u\|_{T}^2 \le  C_3 \int_T u^2. 
\end{equation}
Combining (\ref{stability})  and (\ref{stability2}) together, we have 
\begin{equation}
\label{upperb}
 \|\nabla u\|_{T}^2 \le  C_3 |T|^{-d} \int_T u^2\le C_3 |T|^{-d} V_T \|\bfc\|_2^2 = C_4  \|\bfc\|_2^2 
\end{equation} 
for a positive constant $C_4$ and
\begin{equation}
\label{upperl}
 \|\nabla u\|_{T}^2 \ge  C_1 |T|^{-d} \int_T u^2 \ge C_1/K |T|^{-d}V_T \|\bfc\|_2^2 = C_0  \|\bfc\|_2^2
\end{equation}
for another positive constant $C_0$. 
Then summing over all simplexes in $\triangle$ gives \ref{Wriesz} with new constants $C_1$ and 
$C_2$ with a consideration of
the redundant indices over the common facets of neighboring simplexes of $\triangle$. 
Indeed, it is easy to see that  the largest number of redundant indices is the largest number 
of simplexes sharing a common vertex
$\bfv\in \triangle$. Let $\bfv\in \triangle$ be the vertex sharing by $k_\bfv$ number 
of simplexes in $\triangle$. Let 
$N_\bfv$ be the ball centered at $\bfv$ with radius $|\triangle|/2$. It is easy to see 
all simplexes sharing $\bfv$ are
inside $N_\bfv$. The volume of $N_\bfv$ is less or equal to $A\pi (|\triangle|/2)^d$ for a 
positive constant $A$.  As each of these simplexes contains
the inscribed ball with radius $\rho_\triangle$ which has a volume $A\pi (\rho_\triangle)^d $. 
It follows that number $k_\bfv$
of simplexes in $\triangle$ sharing $\bfv$ is estimated by 
$$
k_v \le \frac{ A\pi (|\triangle|/2)^d}{A\pi \rho_\triangle^d} = \frac{1}{2^d} \beta^d <\infty. 
$$
That is, $C_2= C_4 \frac{1}{2^d} \beta^d$ and $C_1= C_0$.  
\end{proof}

With the estimates (\ref{Wriesz}) in the bivariate setting, Lai and Schumaker in \cite{LS09} 
proved the following result.
\begin{lemma}[Lai and Schumaker, 2009 \cite{LS09}]
\label{thm:LS09}
Let $\omega$ be a cluster of triangles in $\tri$, and let $T \in \omega$ be a triangle, e.g. 
$\omega=T$. 
Then there exists constants  $0 < \sigma < 1$  and
$C$ depending only on  the ratio $C_2/C_1$ in (\ref{Wriesz}) such that if $g$ is a function in  
$\mathcal{W}$ with
\begin{equation}
\langle g, \, w\rangle_{\mathcal{W}}= 0, \qquad \hbox{ for all } w 
\in \mathcal{W} \hbox{ with supp(}w) \subseteq \hbox{star}^k(\omega), \label{new}
\end{equation}
for some fixed $k \ge 1$, then
\begin{equation}
 \|g\cdot \chi_T\|_{\mathcal{W}}  \le C \sigma^k \|g\|_{\mathcal{W}}, 
\label{bnd}
\end{equation}
where $\chi_T$ is the characteristic function of $T$. 
\end{lemma}

Now let us translate the result in Lemma~\ref{thm:LS09} in the setting of $\Omega\subset 
\mathbb{R}^d, d\ge 2$ to have 
\begin{theorem}
\label{main1}
Let $\triangle$ be a simplicial partition of $\Omega\subset \mathbb{R}^d, d\ge 2$. 
Suppose that $\mathcal{W}$ has an $1$-stable basis $\{\mathcal{B}_\xi\}_{\xi\in\mathcal{M}}$ satisfying 
(\ref{Wriesz}). 
Let $\omega$ be a cluster of simplexes in $\tri$ or a simplex, e.g. $\omega=T$. 
Then there exists constants  $0 < \sigma < 1$  and
$C$ depending only on  the ratio $C_2/C_1$ in (\ref{Wriesz}) such that if $g$ is a function 
in  $\mathcal{W}$ with
\begin{equation}
\langle g, \, w\rangle_{\mathcal{W}}= 0, \qquad \hbox{ for all } w 
\in \mathcal{W} \hbox{ with supp(w)} \subseteq \hbox{star}^k(\omega), \label{new}
\end{equation}
for some fixed $k \ge 1$, then
\begin{equation}
 \|g\cdot \chi_T\|_{\mathcal{W}}  \le C \sigma^k \|g\|_{\mathcal{W}}, 
\label{bnd}
\end{equation}
where $\chi_T$ is the characteristic function of $T$. 
\end{theorem}
\begin{proof}
We use the  ideas in the proof of  Lemma 4.2 in \cite{LS09} 
to establish this similar result in Theorem~\ref{main1}. The detail is omitted here. \end{proof} 

Note that the proof of Theorem~\ref{main1} 
is based on the ideas in \cite{GLS02} and \cite{GS02} and an elementary inequality in 
\cite{deBoor76} which is included below.  
\begin{lemma}[de Boor, 1996\cite{deBoor76}] 
\label{deBoor}
If  the sequence $a_0, a_1, . . .$ satisfies 
\begin{equation}
\label{deBoor1}
|a_m|\ge c \sum_{j\ge m}|a_j|, \quad m=0, 1, 2, \cdots, 
\end{equation}
for some $c\in (0, 1)$, then $\lambda = 1- c\in (0, 1)$ and 
\begin{equation}
\label{deBoor2} 
|a_m|\le |a_0| \lambda^m /c, \quad m=0, 1, 2, \cdots 
\end{equation}
\end{lemma}
This equality was called discrete Gronwall inequality as in \cite{DDW75}. Please notice that 
this exponential decay includes a linear decay as a special case.  
For example, let $k\ge 10$ be an integer and let 
$$
a_0=k, a_1=k-1, \cdots, a_k=0, a_{k+j}=0, \forall j\ge 1
$$  
be a sequence. 
Then it is clear that this sequence   decays to zero in a linear fashion, but satisfies the 
de Boor condition (\ref{deBoor1}) for $c=1/k$. 
Indeed, if $m \ge k$, we have (\ref{deBoor1}) as 
the both sides are zero and if $0\le m<k$, we have 
$$
c \sum_{j\ge m}|a_j| =\frac{1}{k} \sum_{i=1}^{k-m} i = \frac{k-m}{k}\frac{k-m+1}{2} \le k-m = a_m. 
$$ 
Lemma~\ref{deBoor} shows that this sequence is of exponential decay in the sense of 
(\ref{deBoor2}). We shall use 
this sequence to explain the decay property of the GBC functions over a rectangular domain. 
That is, the GBC functions over a rectangular domain linearly decay to zero and 
satisfy the decay property (\ref{deBoor2}).

We now explain how to use Theorem~\ref{main1} for establishing the e-locality of 
our numerical harmonic GBC functions $S_i$ and $R_i$. 
First of all, we explain the e-locality of $S_i= S_{i,0}+G_i$. We start with $G_i$.   
Let us fix a simplicial partition $\triangle_k$ which is a  refinement of $\triangle$ and 
fix a spline space, say $S^1_n(\triangle_k)$ with $n\ge 5$ in $\mathbb{R}^2$ or $n\ge 9$ 
in $\mathbb{R}^3$ and etc. 
Then we can construct $G_i\in S^1_n(\triangle)$ satisfying the boundary condition, i.e.
$G_i|_{\partial \Omega}= \ell_i$ as follows.   
Without loss of generality, we may assume that $\Omega$ is a star-shaped polygon in 
$\mathbb{R}^2$. 
There is a  center $\bfv_c$ which can be connected to all vertices of $\Omega$. Let us form 
an initial 
triangulation in this way, say $\triangle_0$. Then our $\triangle_k$ is a uniformly refined 
triangulation of $\triangle_0$. Let us say $\triangle$ is the third refinement of $\triangle_0$. 

If we use a continuous spline space $S^0_n(\triangle_k)$ with $n\ge 1$,   
we can easily construct $G_i$ at boundary vertex $\bfv_i$. Indeed, let us use 
Figure~\ref{fig0} to show  
spline coefficients of $G_i$ in $S^1_5(\triangle_0)$ over the domain points on three 
triangles. The remaining  
spline coefficients of $G_i$ over other triangles of $\triangle_0$ are all zero.  

\begin{figure}[thbp]
\begin{picture}(500,220)(-250,0)
\unitlength=0.6pt
\thicklines
\multiput(-190,50)(50,0){5}{\line(1,0){30}}
\multiput(-195,60)(20,40){5}{\line(1,2){10}}
\multiput(41,62)(-30,40){5}{\line(-3,4){15}}
\put(52,60){\line(1,5){48}}
\put(-92,252){\line(4,1){192}}
\put(-202,62){\line(-1,6){48}}
\put(-106,254){\line(-3,2){144}}
\put(-205,45){$1$} 
\put(-155,45){$\frac{4}{5}$}
\put(-105,45){$\frac{3}{5}$}

\put(-55,45){$\frac{2}{5}$}
\put(-5,45){$\frac{1}{5}$}
\put(55,45){$0$}
\put(-185,85){$0$}
\put(-135,85){$0$}

\put(-85,85){$0$}
\put(-35,85){$0$}
\put(15,85){$0$}

\put(-220,107){$\frac{4}{5}$}
\put(-193,147){$0$}

\put(-173,187){$0$}
\multiput(-153,227)(20,40){2}{$0$}
\put(-230,167){$\frac{3}{5}$}
\multiput(-203,207)(20,40){3}{$0$}
\put(-240,227){$\frac{2}{5}$}
\put(-213,267){$v$}
\put(-193,307){$0$}
\multiput(-250,287)(20,40){1}{$\frac{1}{5}$}
\put(-230,327){$0$}

\put(-253,347){$0$}

\multiput(57,97)(-30,40){2}{$0$}
\put(-13,177){$0$}
\multiput(-43,217)(-30,40){2}{$0$}
\multiput(67,147)(-30,40){4}{$0$}
\multiput(77,197)(-30,40){3}{$0$}
\multiput(87,247)(-30,40){2}{$0$}
\put(97,297){$0$}
\put(-165,125){$0$}
\put(-15,125){$0$}
\put(-115,128){$0$}
\put(-65,128){$0$}
\put(-145,165){$0$}
\put(-95,165){$0$}
\put(-45,165){$0$}
\put(-125,205){$0$}
\put(-75,205){$0$}
\put(-105,245){$0$}
\put(-220,25){${\bf v}_{1}$}
\put(50,25){${\bf v}_{2}$}
\put(-105,260){${\bf v}_c$}
\put(-270,360){${\bf v}_{n}$}
\end{picture}
\caption{A Construction of $G_1$ over the two triangles sharing $\bfv_1$ using continuous spline space 
$S^0_5(\triangle)$  \label{fig0}}
\end{figure}

Next let us explain how to construct $G_i\in S^1_5(\triangle_0)$. 
 Once we have $G_i$, we can rewrite $G_i$ over $S^1_5(\triangle_k)$ 
due to the nestedness of our spline spaces.    
On $S^1_5(\triangle_0)$, $G_1$ can be constructed by specifying the coefficients as shown in Figure~\ref{triangle} 
over domain points of degree $5$ over $\triangle_0$.    
Note that $\bfv_1, \bfv_2, \cdots, \bfv_N$ are boundary vertices. $\bfv_c$ is an interior vertex. 
$u, v, X, Y\ge 0$ are coefficients which are dependent on $C^1$ smoothness conditions.  
That is, we use $C^1$ smoothness condition 
connected to the coefficients $1$, two $4/5$ and $Y$ to find $Y$ first if the edge 
$\bfe_1=\langle\bfv_1, \bfv_2\rangle$ and the edge $\langle \bfv_1, \bfv_{N}\rangle$ are not parallel.  
Then we find two positive values $X$ by solving the $C^1$ smoothness condition connecting the coefficients $Y$, $0$ 
and $X, X$. If the edge 
$\bfe_1=\langle\bfv_1, \bfv_2\rangle$ and the edge $\langle \bfv_1, \bfv_{N}\rangle$ are parallel,  
we simply choose $Y=0$ and hence, $X=0$.  
Next we use the $C^1$ smoothness condition to find $u$ and then $v$.  
The remaining coefficients of $G_1$ are all zero which we do not show
in Figure~\ref{triangle}.  Hence, this function $G_1$ is in $S^1_5(\triangle_0)$.  

\begin{figure}[thbp]
\begin{picture}(500,220)(-250,0)
\unitlength=0.6pt
\thicklines
\multiput(-190,50)(50,0){5}{\line(1,0){30}}
\multiput(-195,60)(20,40){5}{\line(1,2){10}}
\multiput(41,62)(-30,40){5}{\line(-3,4){15}}
\put(52,60){\line(1,5){48}}
\put(-92,252){\line(4,1){192}}
\put(-202,62){\line(-1,6){48}}
\put(-106,254){\line(-3,2){144}}
\put(-205,45){$1$} 
\put(-155,45){$\frac{4}{5}$}
\put(-105,45){$\frac{3}{5}$}

\put(-55,45){$\frac{2}{5}$}
\put(-5,45){$\frac{1}{5}$}
\put(55,45){$0$}
\put(-185,85){$Y$}
\put(-135,85){$X$}

\put(-85,85){$0$}
\put(-35,85){$u$}
\put(15,85){$\frac{1}{5}$}

\put(-220,107){$\frac{4}{5}$}
\put(-193,147){$X$}

\put(-173,187){$0$}
\multiput(-153,227)(20,40){2}{$0$}
\put(-230,167){$\frac{3}{5}$}
\multiput(-203,207)(20,40){3}{$0$}
\put(-240,227){$\frac{2}{5}$}
\put(-213,267){$v$}
\put(-193,307){$0$}
\multiput(-250,287)(20,40){2}{$\frac{1}{5}$}

\put(-253,347){$0$}

\multiput(57,97)(-30,40){2}{$0$}
\put(-13,177){$0$}
\multiput(-43,217)(-30,40){2}{$0$}
\multiput(67,147)(-30,40){4}{$0$}
\multiput(77,197)(-30,40){3}{$0$}
\multiput(87,247)(-30,40){2}{$0$}
\put(97,297){$0$}
\put(-165,125){$0$}
\put(-15,125){$0$}
\put(-115,128){$0$}
\put(-65,128){$0$}
\put(-145,165){$0$}
\put(-95,165){$0$}
\put(-45,165){$0$}
\put(-125,205){$0$}
\put(-75,205){$0$}
\put(-105,245){$0$}
\put(-220,25){${\bf v}_{1}$}
\put(50,25){${\bf v}_{2}$}
\put(-105,260){${\bf v}_c$}
\put(-270,360){${\bf v}_{n}$}
\end{picture}
\caption{A Construction of $G_1\in S^1_5(\triangle_0)$ over three triangles sharing $\bfv_1$  \label{triangle}}
\end{figure}

In the same way, we can construct $G_i\in S^0_n(\triangle_k)$ or in $S^1_n(\triangle_k)$  
for $\triangle_k$ in the higher dimensional setting. 
We leave the detail to the interested reader.  It is clear
that $G_i$ is of exponential decay.  To show the e-locality of $S_i$, 
we only need to discuss the e-locality of
$S_{i,0}$ with $S_{i,0}=S_i- G_i$.  

\begin{figure}[htbp]
\centering
\includegraphics[width = 0.5\textwidth]{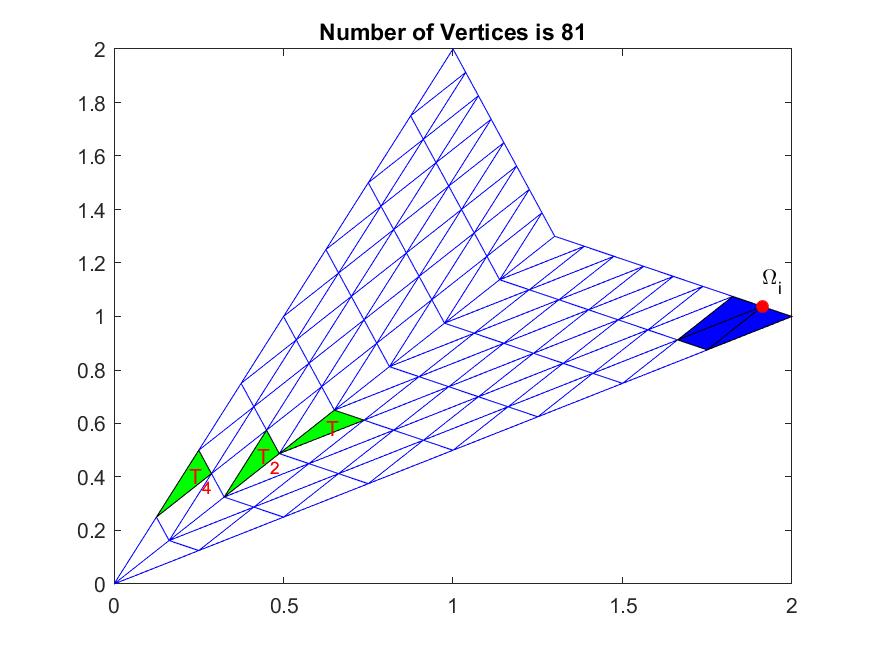}
\caption{Illustration of $T$ whose $\hbox{star}^5(T)\subset \Omega_i^c$ and 
$T_2$ and $T_4$. Note that $T_3$ is not shown. It can be any triangle between $T_2$ and 
$T_4$ as long as $T_3\cap T_2\not=\emptyset$
and $T_3\cap T_4\not=\emptyset$.    \label{Ltriangle}}
\end{figure}

To this end,  let $\Omega_i=\hbox{star}^1(\bfv_i)$ be the union of all simplexes 
which are connected to the boundary vertex $\bfv_i$. According to the construction above, the 
support of $G_i$ is contained in $\Omega_i$ in both $C^0$ and $C^1$ cases.   We say
a simplex $T$ is L-simplex away from the boundary $\partial \Omega$ if $L$ is the 
smallest integer such that 
there exist  $L$ simplexes $T_1, \cdots, T_{L+1}$ in $\triangle$ such that $
T=T_1$ and $\bar{T}_{L+1}\cap \partial \Omega\not=\emptyset$ as well as  $\bar{T}_j\cap
\bar{T}_{j+1}\not=\emptyset$ for $j=1, \cdots, L$. Furthermore, we say $T$ is $L$-simplex 
away from the boundary $\partial \Omega$ in the opposite direction of $\bfv_i$ if there 
exists an integer $k$ such 
that $\hbox{star}^k(T)\subset \Omega_i^c$, the complement of $\Omega_i$ 
in $\Omega$ and $T_i\in \hbox{star}^{k+\ell_i}(T)\subset \Omega_i^c$ for $\ell_{i+1}\ge 
\ell_{i}\ge 0$ for $i=1, \cdots, L+1$. See Figure~\ref{Ltriangle} for $\bfv_i$ (the red dot),  
$\Omega_i$, $T=T_1$,  $T_2$ and $T_4$, 
where $T_3$ can be any triangle between $T_2$ and $T_4$ as long as $T_3\cap T_2\not=\emptyset$
and $T_3\cap T_4\not=\emptyset$. There are three choices for $T_3$. 

We are now ready to state and prove one of the main results in this paper. 

\begin{theorem}
\label{main}
Fix a boundary vertex $\bfv_i\in \Omega$. Suppose that $\Omega_i=\hbox{star}(\bfv_i)\not=\Omega$.
Let  $\bfv\in \Omega$ be a point in $\Omega\setminus\Omega_i$.  
Let $T\in \triangle$ contain $\bfv$, i.e. $\bfv\in T$. 
Suppose $k\ge 1$ is the integer such that $\hbox{star}^k(T) \subset \Omega_i^c$, 
the complement of $\Omega_i$ 
in $\Omega$ and $L\ge 1$ is another integer such that 
$T$ is $L$ simplexes away from the partial boundary $\partial \Omega\backslash \Omega_i$ in the
opposite direction of $\bfv_i$. Then there exists  a positive
constant $K>0$ independent of $L$ and $\sigma\in (0, 1)$ such that 
\begin{equation}
\label{eq:main2}
|S_i(\bfv)|\le K (L+1) \sigma^k. 
\end{equation}
\end{theorem}
\begin{proof}
It is easy to see from the construction of $G_i$ above that $G_i$ is only supported over 
$\Omega_i$. Since $\Omega_i\subset \Omega$ and $\Omega_i \not=\Omega$, 
$G_i$ has e-locality. That is, $G_i$ decays to zero outside of $\Omega_i$.  
We only need to show that $S_{i,0}=S_i-G_i$ has e-locality. 
According the assumptions, $\hbox{star}^k(T) \cap \Omega_i =\emptyset$. 
By (\ref{weak}), we use $S_{i,0}$ for $g$ in Theorem~\ref{main1} 
and use (\ref{weak}) to have  
$$
\langle g, \, w\rangle_{\mathcal{W}}:=\langle \nabla S_{i,0}, \nabla w\rangle = -\langle \nabla G_i, 
\nabla w\rangle = 0,  
$$
for all $ w \in \mathcal{W}$ with $\hbox{sup}(w) \subseteq \hbox{star}^k(T)$.  
Therefore, we can use the result from Theorem~\ref{main1} above to conclude that there exists a 
positive number $\sigma\in (0, 1)$ such that 
\begin{equation}
\label{main2}
\|\nabla S_{i,0}\cdot \chi_T\|\le C \sigma^k \|\nabla S_{i,0}\|. 
\end{equation}

We now show $\|S_{i,0}\|_\infty$ can be bounded by the left-hand side of (\ref{main2}).  
To do so, let us use induction on $L\ge 1$. If $L=1$, then $T$ is one simplex away 
from the partial boundary $\partial \Omega \backslash \Omega_i$ in the opposite direction of $\bfv_i$. 
Let  $\bfv \in T\cap T_2$ and $\bfu\in T_2\cap (\partial \Omega\backslash \Omega_i)$. 
That is, $\bfu$ on the partial boundary  
and hence $S_{i,0}(\bfu)=0$.  By Taylor expansion with remainder, we have
$$
0= S_{i,0}(\bfu)= S_{i,0}(\bfv) + \nabla S_{i,0}(\bfw)(\bfu- \bfv)
$$
for a point $\bfw\in T$ in the line segment between $\bfu$ and $\bfv$. It follows that 
$$
|S_{i,0}(\bfv)|=|\nabla S_{i,0}(\bfw)(\bfu- \bfv)| \le |\triangle| \| \nabla S_{i,0}\|_{\infty, T_2} \le C\|\nabla S_{i,0}\|_{L^2(T_2)} 
=C\|\nabla S_{i,0}\cdot \chi_{T_2}\|
$$
by Theorem~1.1 of \cite{LS07} over triangle $T_2$. Similarly, we have the same estimate for polynomials over simplex
$T$ in $\mathbb{R}^d$ if $d>2$. We now use (\ref{main2}) to conclude 
$$
|S_{i,0}(\bfv)| \le C \sigma^{k} \|\nabla S_{i,0}\|. 
$$ 
As $\|\nabla S_{i,0}\|=\|\nabla (S_i-G_i)\|$ is very close to $\|\nabla (\phi_i-G_i)\|$ by (\ref{cea}), 
there exists a constant $K>0$ dependent on $\nabla \phi_i$ such that 
$$
|S_{i,0}(\bfv)| \le C \sigma^k \|\nabla S_{i,0}\| \le C \sigma^k \|\nabla (S_i-\phi_i)\| + 
C\sigma^k \|\nabla (\phi_i-G_i)\|  \le K \sigma^k. 
$$
Hence, $|S_i(\bfv)|=  \|S_{i,0}(\bfv)| \le K \sigma^k. $ 

Next assume that when $T$ is $(L-1)$-simplex away from the partial boundary of $\Omega$ in the 
opposite direction of
$\bfv_i$, we have the desired estimate. Let us
consider the case when $T$ is $L$-simplex away from the partial boundary of $\Omega$ 
in the opposite direction of $\bfv_i$. 
Let $\bfu_L\in T_{L}$ be the intersection of $T_{L}$ and $T_{L+1}$ and $T_{L+1}$ intersects the 
partial boundary. By the argument above, we have
$$
|S_{i,0}(\bfu_L)| \le K \sigma^{k+\ell_L}. 
$$ 
for an integer $\ell_L$ with $\ell_L\ge 0$. 
For $\bfu_j\in T_j\cap T_{j+1}$ with $j=1,\cdots, L-1$ with $\bfu_1=\bfv$,  
we use Taylor expansion again to have
$$
S_{i,0}(\bfu_j)= S_{i,0}(\bfu_{j+1}) + \nabla S_{i,0}(\bfw)(\bfu_j- \bfu_{j+1})
$$
for an appropriate $\bfw\in T_j$. It follows that 
$$
|S_{i,0}(\bfu_j)|\le |S_{i,0}(\bfu_{j+1})|+ |\tri| \|S_{i,0}\|_{\infty, T_j} 
\le |S_{i,0}(\bfu_{j+1})| + C\|\nabla S_{i,0}\|_{L^2(T_j)}
$$
by Theorem~1.1 of \cite{LS07} again and 
$$
|S_{i,0}(\bfv)| \le \cdots \le K \sigma^{k+\ell_L} + \sum_{j=1}^L C\|\nabla S_{i,0}\|_{L^2(T_j)}.
$$ 
As $T_j$ is in the opposite direction of $\bfv_i$, similar to 
(\ref{main2}), we have $\|\nabla S_{i,0}\|_{L^2(T_j)}\le C\sigma^{k+\ell_j} \|\nabla S_{i,0}\|$. 
We put these inequalities together to have the desired
result follows. We have thus completed the proof. 
\end{proof}

When $S_i$ is a GBC function over a rectangular domain $\Omega$, the above arguments show that $S_i$ satisfies 
the decay problem (\ref{eq:main2}) although $S_i$ decays linearly. We can say that 
$S_i$ decays in the sense of  (\ref{deBoor2}). Except for this pathological example, the GBC 
functions have  indeed e-locality as shown in Example~\ref{ex1}, i.e. Figures~\ref{quad2} and \ref{quad3} and 
numerical experiments the researchers and practitioners have already observed.

\begin{example}
\label{ex1}
Consider a quadrilateral and a triangulation $\triangle$ shown on Figure~\ref{quad} 
\begin{figure}[htbp]
\centering
\includegraphics[width = 0.4\textwidth]{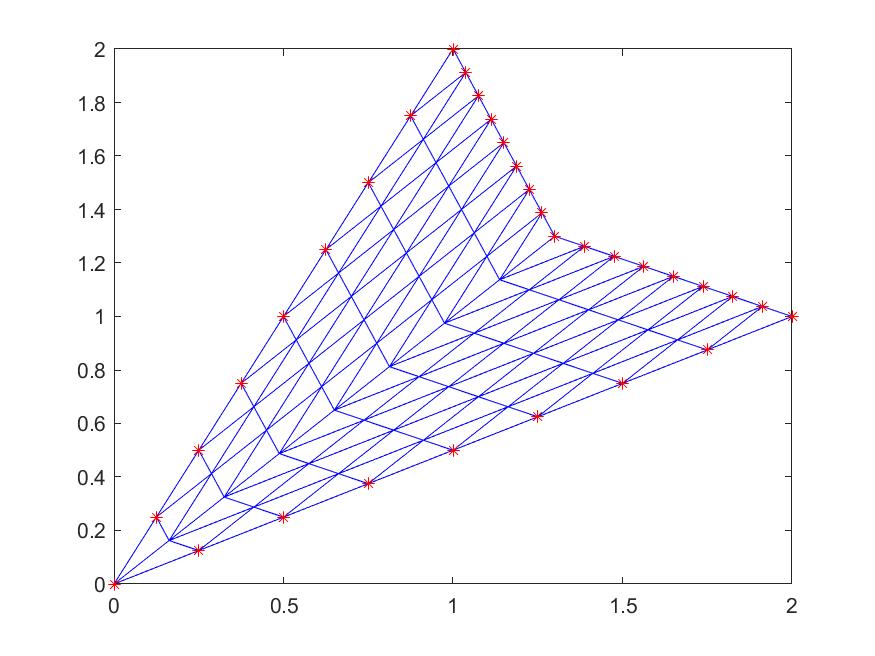}
\caption{A triangulation of polygon $\Omega$ with boundary vertices \label{quad}}
\end{figure} 

We use $S^1_5(\triangle)$ to approximate harmonic GBC functions. There are 32 GBC functions 
associated with boundary vertices shown in red in Figure~\ref{quad}.  For convenience, we 
show 6 of them to illustrate that these functions clearly have e-decay property. 

\begin{figure}[htbp]
\centering
\includegraphics[width = 0.3\textwidth]{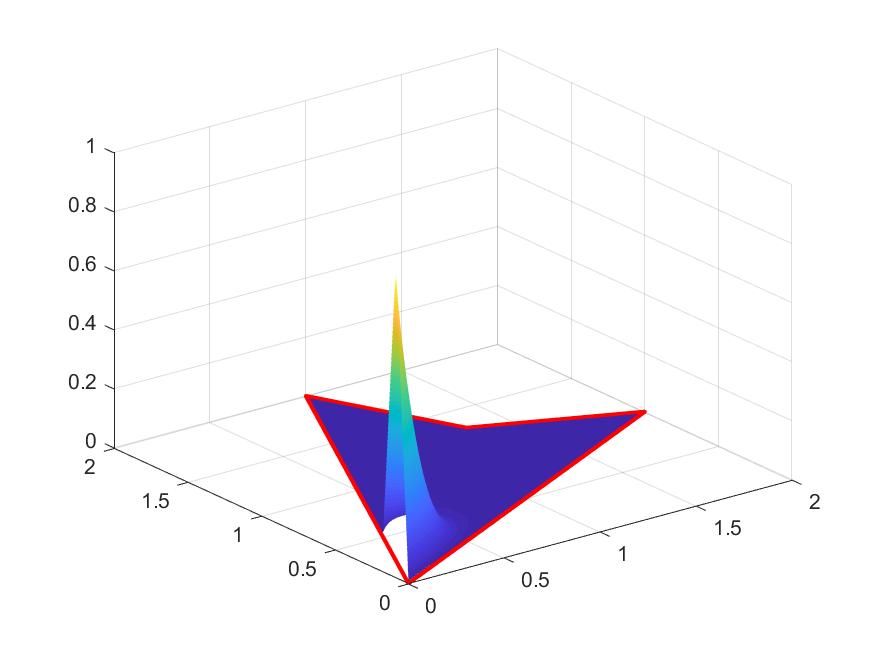}
\includegraphics[width = 0.3\textwidth]{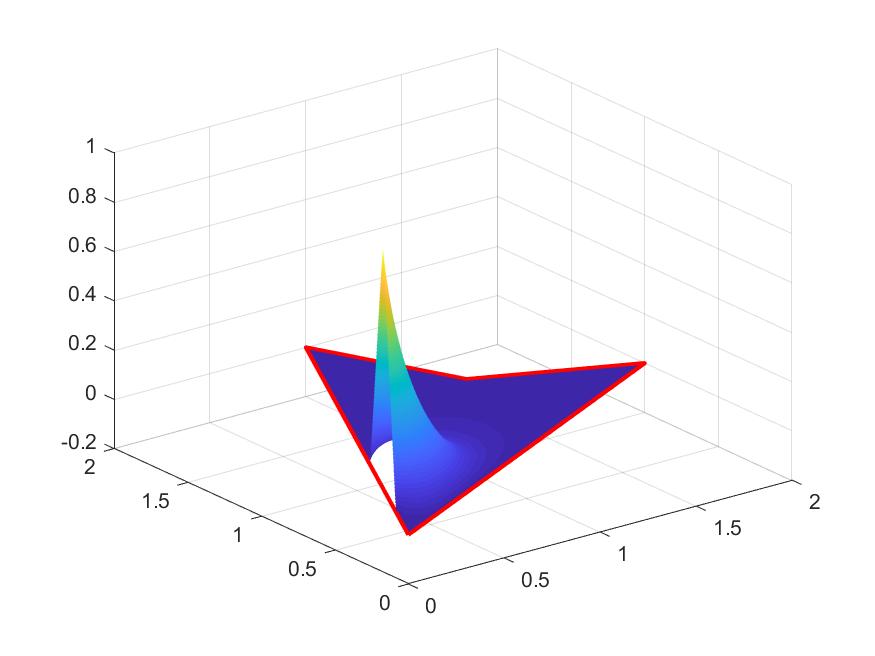}
\includegraphics[width = 0.3\textwidth]{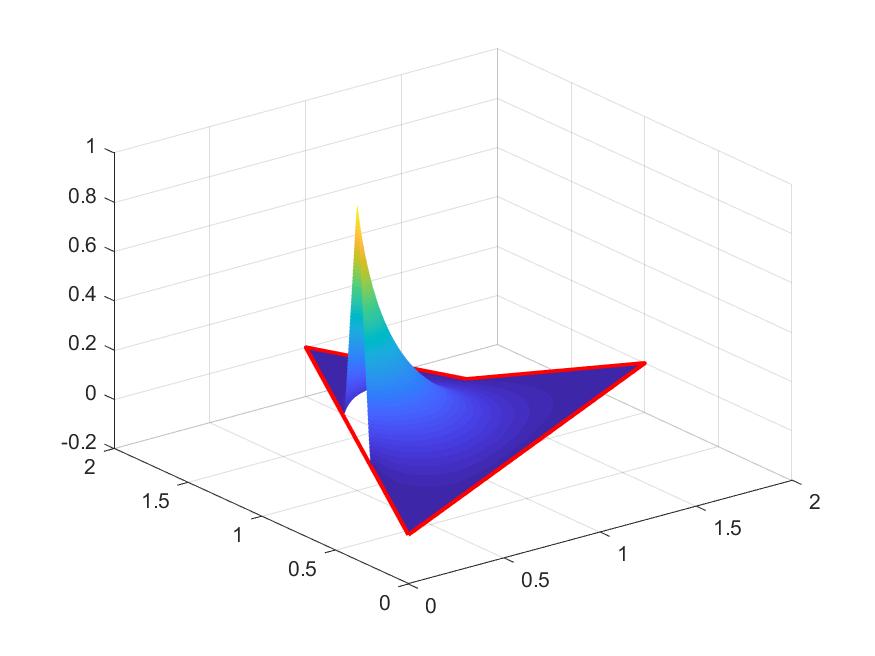}
\caption{Spline Approximation of Various Harmonic GBC Functions on the Boundary of a Polygon. \label{quad2}}
\end{figure} 

\begin{figure}[htbp]
\centering
\includegraphics[width = 0.3\textwidth]{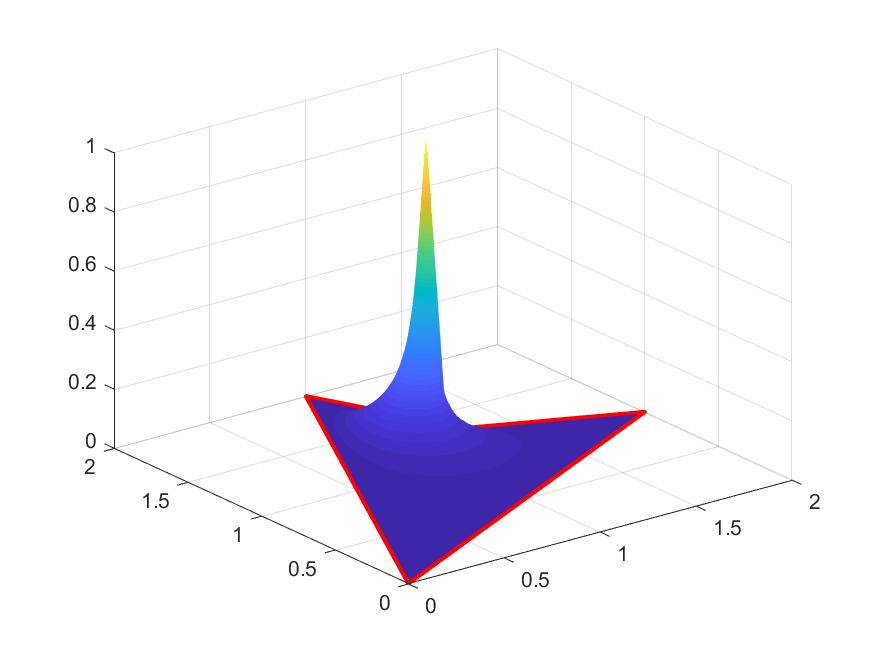}
\includegraphics[width = 0.3\textwidth]{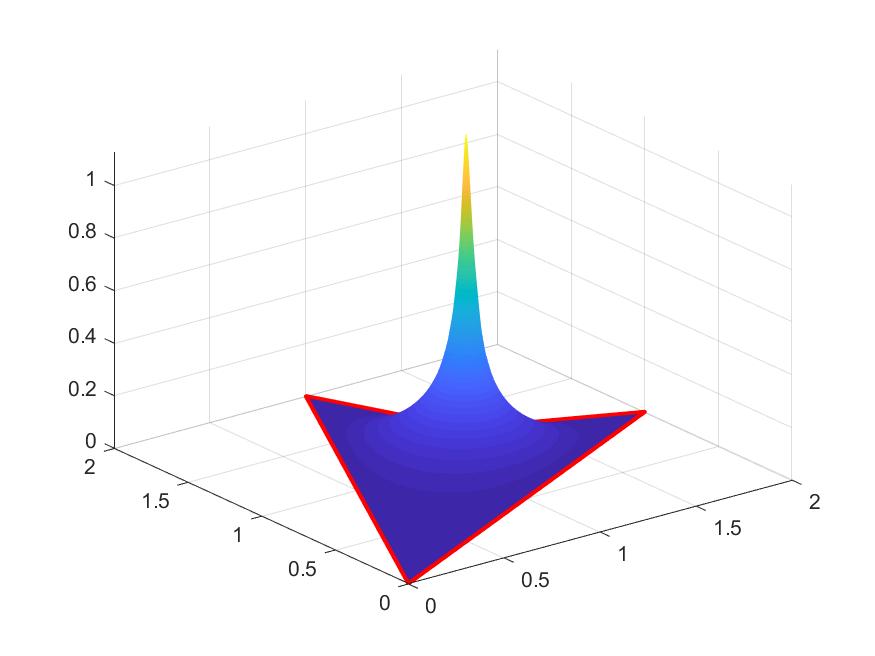}
\includegraphics[width = 0.3\textwidth]{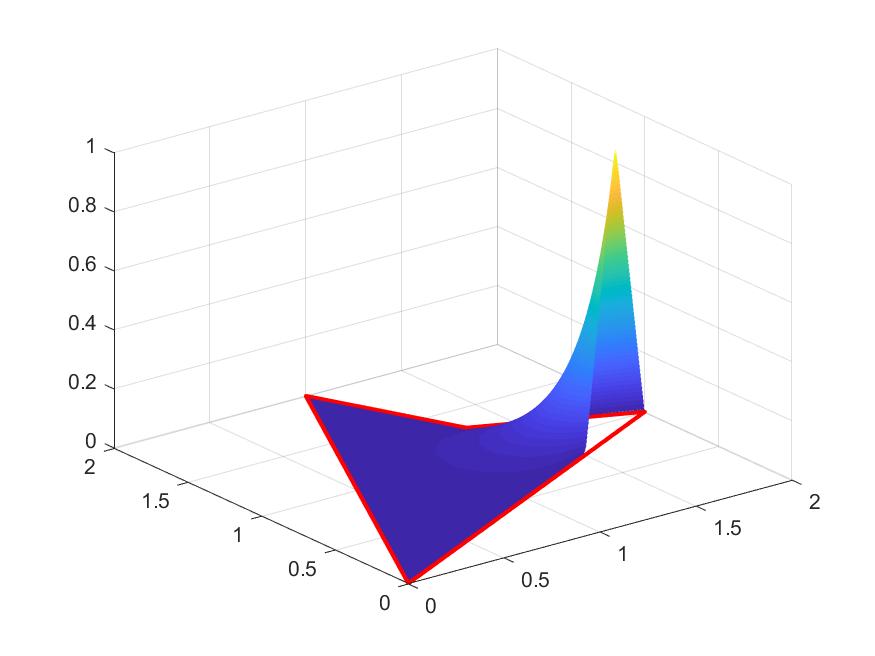}
\caption{Spline Approximation of Various Harmonic GBC Functions on the Boundary of a Polygon \label{quad3}}
\end{figure} 

From graphs above, we clearly see the e-decay property of various spline harmonic GBC 
functions. Similar for other 26 GBC functions $S_i$ which are not shown here. 
\end{example}

Next we show the e-locality of spline harmonic function $R_i$. Let $\omega_i$ be the support of function $h_i$
and $\widetilde{\Omega}_i$ be the union of all $\hbox{star}^1(t),t\in \omega_i$.  Similar to $S_i$, we have  
\begin{theorem}
\label{main3}
Let  $\bfv\in \Omega$ be a point in $\Omega$ which is not in $\widetilde{\Omega}_i$.  
Let $T\in \triangle$ contain $\bfv$, i.e. $\bfv\in T$. 
Suppose that  $k\ge 1$ is an integer  such that $\hbox{star}^k(T) \subset \widetilde{\Omega}_i^c$, the complement of $\widetilde{\Omega}_i$ 
in $\Omega$. Suppose that $T$ is $L$ simplexes away from the  boundary $\partial \Omega$ opposite to $\bfv_i$. 
Then there exists  a positive
constant $K>0$ independent of $L$ and $\sigma\in (0, 1)$ such that 
\begin{equation}
\label{eq:main3}
|R_i(\bfv)|\le K (L+1) \sigma^k. 
\end{equation}
\end{theorem}
\begin{proof}
For any triangle $T\in \Omega \backslash \widetilde{\Omega}_i$, 
we know $\hbox{star}^k(T) \cap \widetilde{\Omega} =\emptyset$. 
By (\ref{weak2}), we use $R_i$ for $g$ in Theorem~\ref{main1} and (\ref{new}) to have 
$$
\langle g, \, w\rangle_{\mathcal{W}}:=\langle \nabla R_i, \nabla w\rangle = \langle h_i, w\rangle = 0,  
$$
for all $ w \in \mathcal{W}$ with $\hbox{sup}(w) \subseteq \hbox{star}^k(T)$.  
Therefore, it follows from Theorem~\ref{main1} that there exists a 
positive number $\sigma\in (0, 1)$ such that 
\begin{equation}
\label{main5}
\|\nabla R_i\cdot \chi_T\|\le C \sigma^k \|\nabla R_i\|. 
\end{equation}
Next we use induction on $L\ge 1$.  The argument is the same as the proof of Theorem~\ref{main}.  We omit the detail 
here to complete the proof.  
\end{proof}

\begin{example}
Consider a quadrilateral and a triangulation $\triangle$ shown on Figure~\ref{quad}.  
We use $S^1_5(\triangle)$ to approximate interior-harmonic GBC functions. There are 49 GBC functions 
associated with interior vertices shown  in Figure~\ref{quad}.  For convenience, we 
show 6 of them in Figure~\ref{quad4} and \ref{quad5}, where these functions also have e-decay property. 

\begin{figure}[htbp]
\centering
\includegraphics[width = 0.3\textwidth]{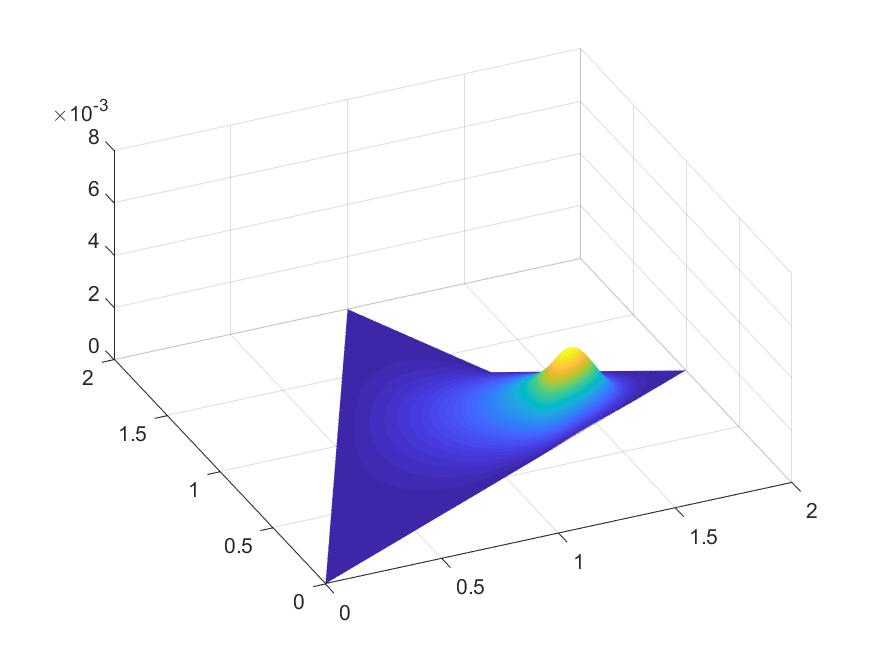}
\includegraphics[width = 0.3\textwidth]{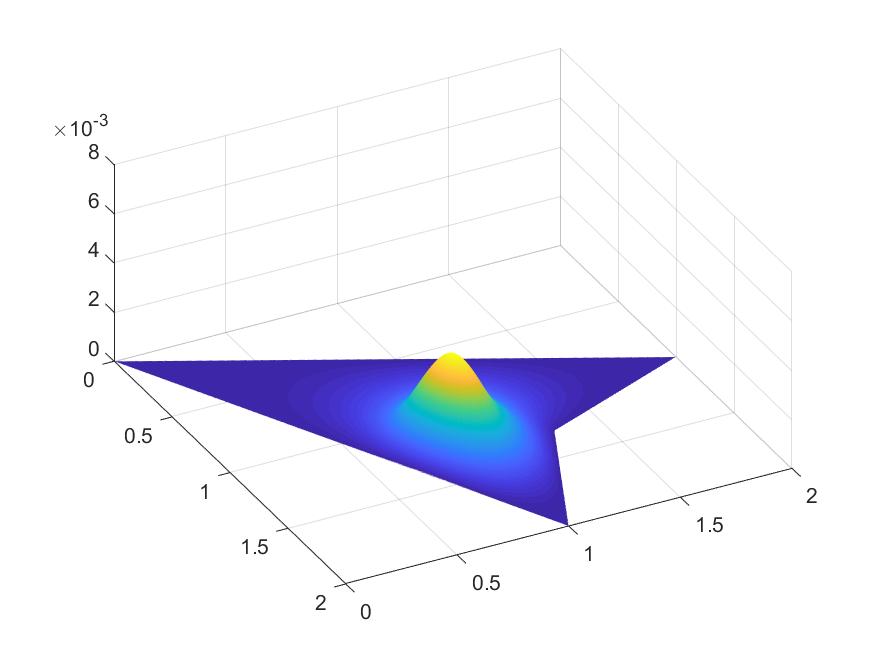}
\includegraphics[width = 0.3\textwidth]{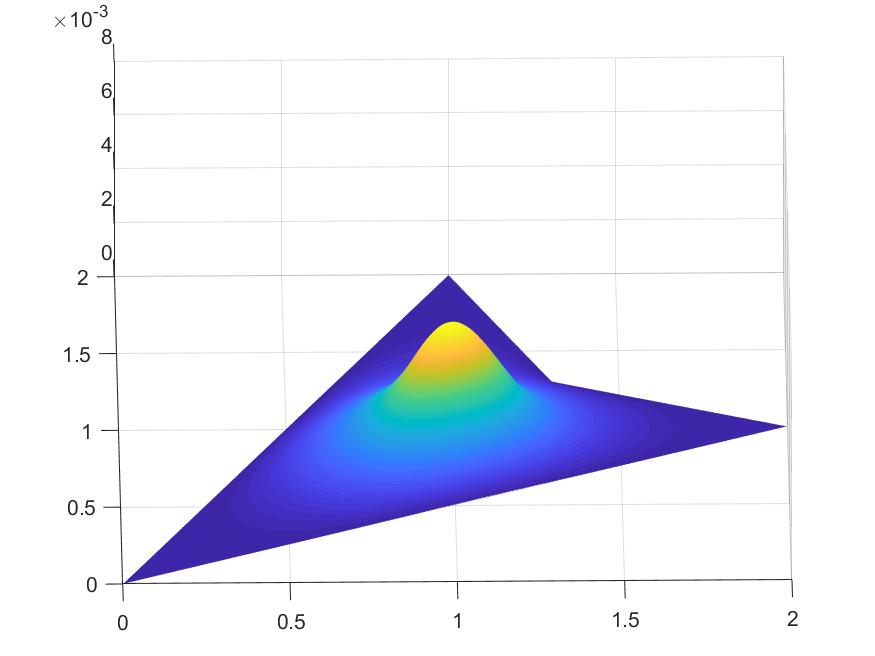}
\caption{Spline Approximation of Various Interior-Harmonic GBC Functions over the Polygon\label{quad4}}
\end{figure} 

\begin{figure}[htbp]
\centering
\includegraphics[width = 0.3\textwidth]{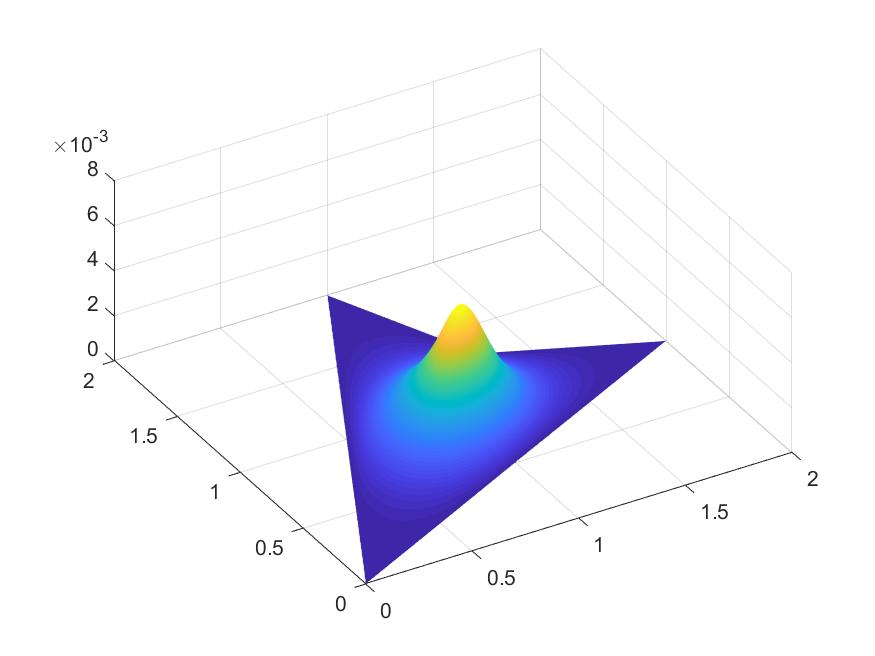}
\includegraphics[width = 0.3\textwidth]{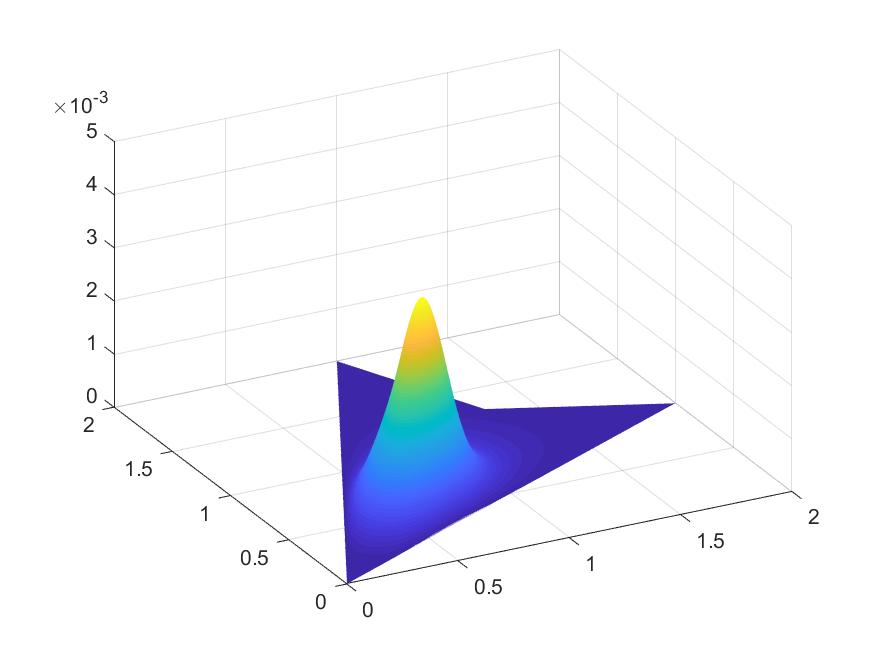}
\includegraphics[width = 0.3\textwidth]{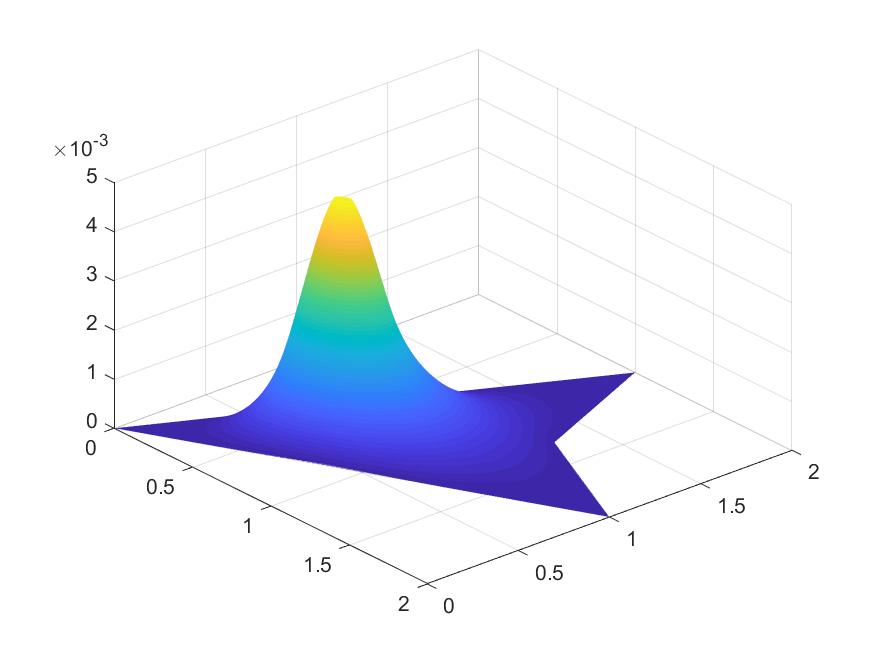}
\caption{Spline Approximation of Various Interior-Harmonic GBC Functions \label{quad5}}
\end{figure} 

From graphs above, we clearly see the e-decay property of various spline interior harmonic GBC 
functions. Similar for remaining 43 GBC functions $R_i$ which are not shown here.
\end{example}

\section{Local Boundary and Interior GBC Functions}
\begin{figure}[htbp]
\centering
\includegraphics[width = 0.4\textwidth]{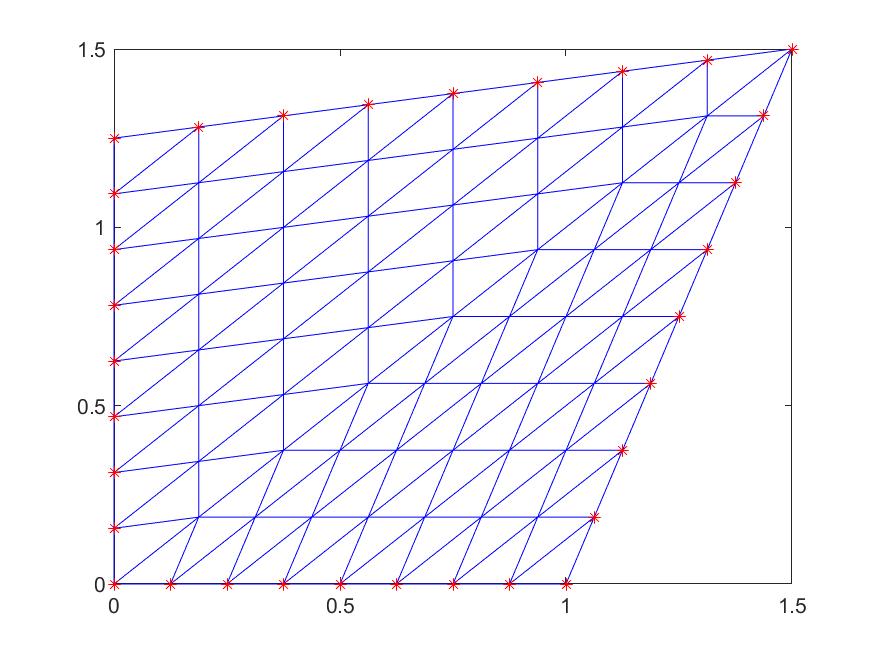}
\includegraphics[width = 0.4\textwidth]{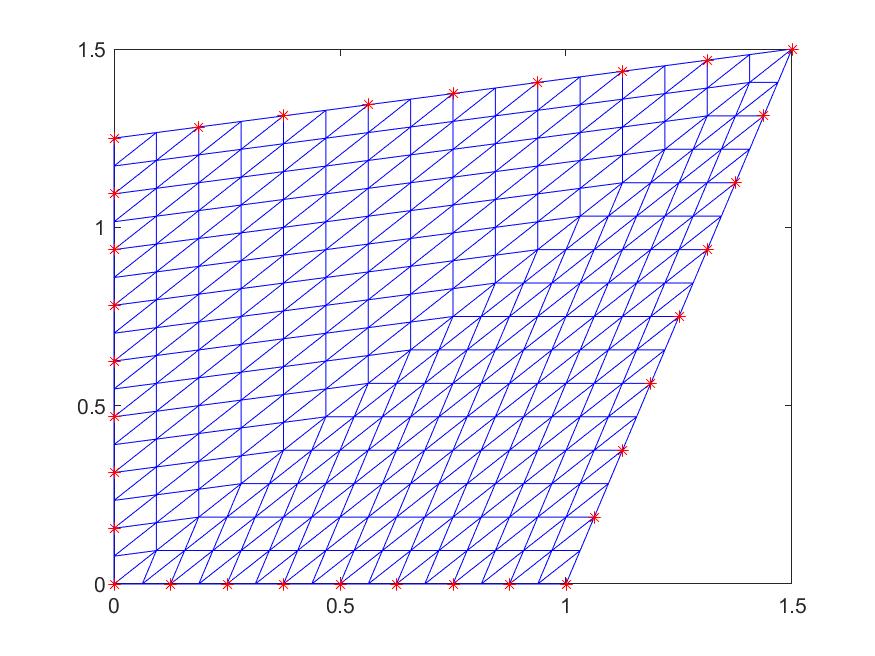}
\caption{A triangulation of polygon $\Omega$ with boundary vertices and its refinement
\label{cquad}}
\end{figure}

\begin{figure}[htbp]
\centering
\includegraphics[width = 0.3\textwidth]{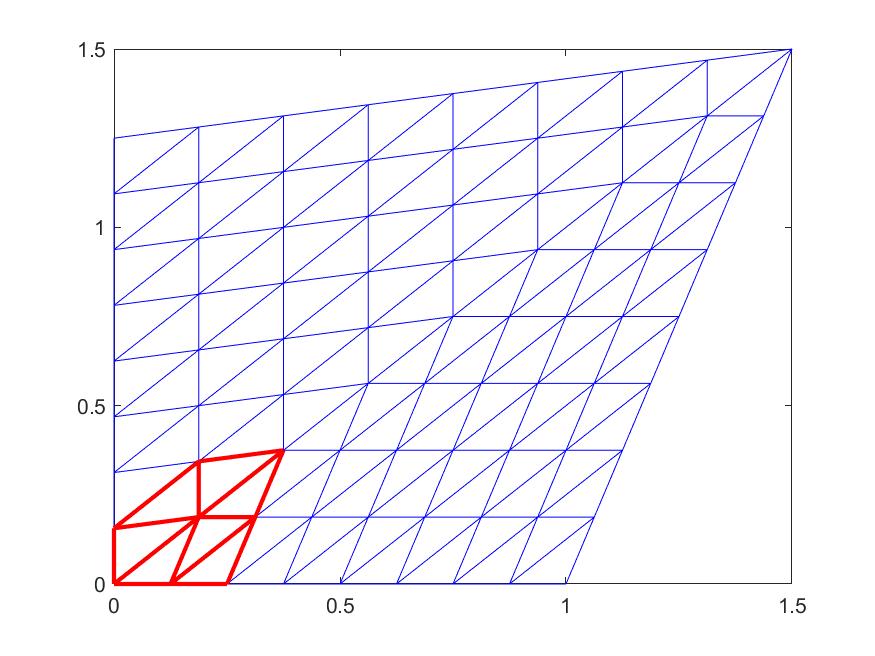}
\includegraphics[width = 0.3\textwidth]{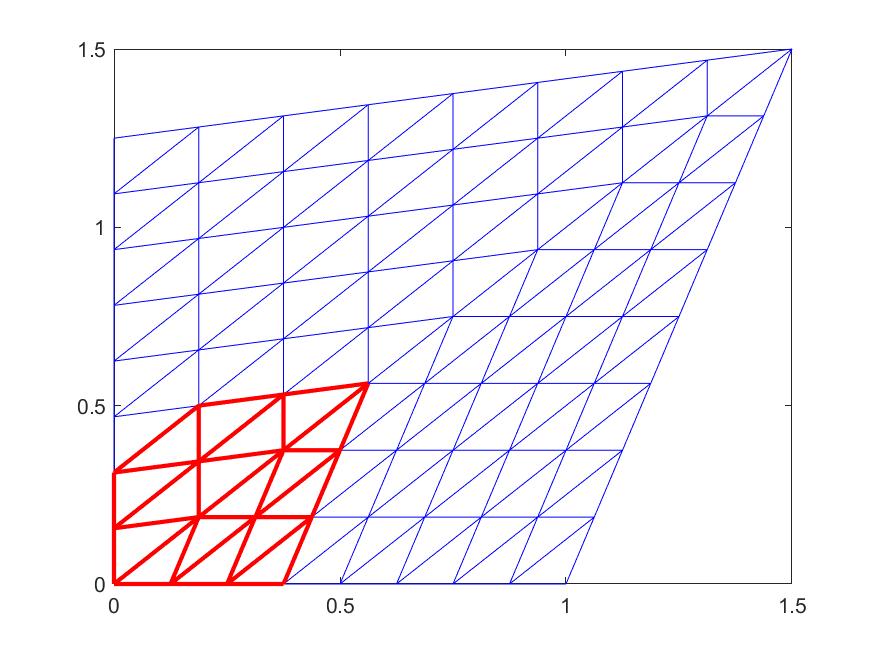}
\includegraphics[width = 0.3\textwidth]{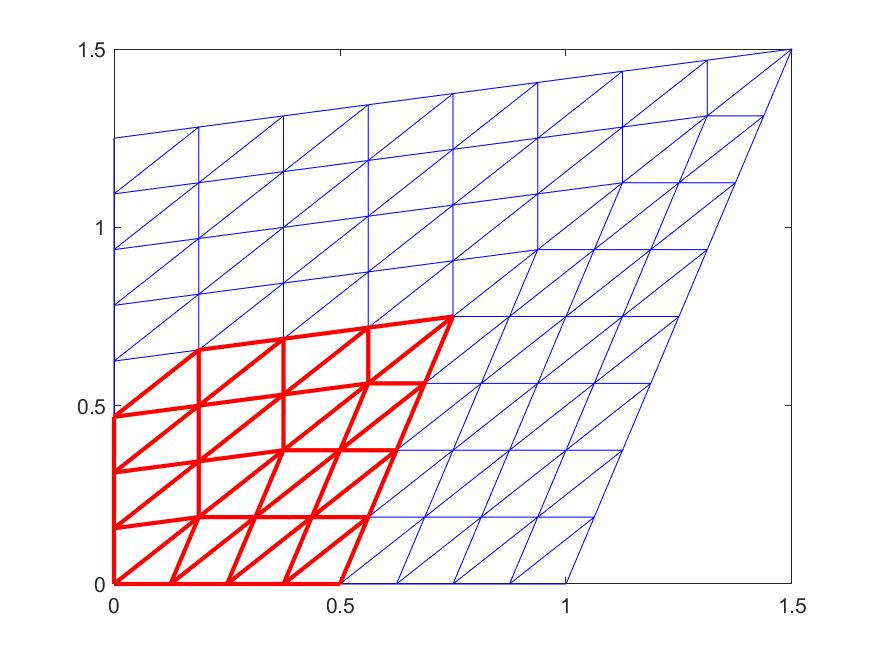}
\caption{Various star$^k(t_0)$ for $k=1,2,3$ \label{rings}}
\end{figure}

In this section, we explore the e-locality of GBC functions using a numerical method. 
Mainly, it is interesting to know
how small $\sigma\in (0,1)$ is. As the $\sigma$ is small, it makes sense to use GBC functions 
supported over star$^k(\bfv)$ to approximate the GBC function with supporting vertex $\bfv$ for 
some fixed $k\ge 1$. 
 For convenience we simply use a convex polygon(as shown on the left in Figure~\ref{cquad}) 
to numerically 
compute a spline approximation of  GBC functions (boundary and interior) globally and 
locally, respectively based on the computational triangulation (the right as shown in 
Figure~\ref{cquad}). That is, for each boundary vertex (the red points 
given in Figure~\ref{cquad}), we compute the standard GBC functions. For each interior
vertex (the vertices inside the polygon on the left graph), we compute the interior
GBC function. Here we say k-local GBC functions if the GBC function with supporting vertex 
$\bfv$ is computed based on the star$^k(\bfv)$ for each $k\ge 1$. In general, we should 
use $k\ge 2$.   
Fix a boundary vertex $\bfv_i$ with $i\in V_B$.  Let $t_i$ be a triangle in $\triangle$ with 
$\bfv_i$ as one of its 
vertex.  For each ring number $k\ge 1$, let $V_k, T_k$ be the triangulation of star$^k(t_i)$ 
shown in Figure~\ref{rings}.

\begin{figure}[htbp]
\centering
\includegraphics[width = 0.3\textwidth]{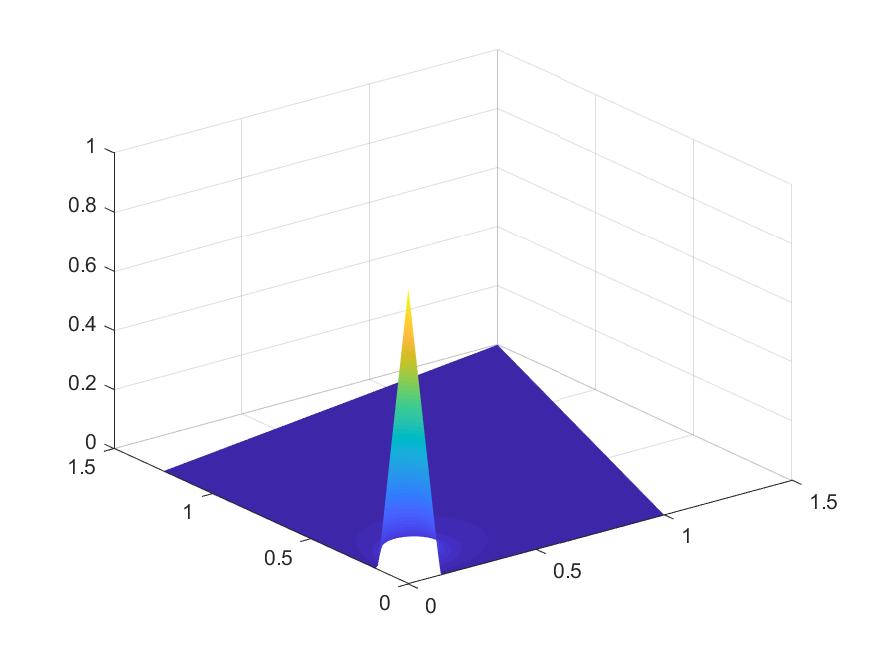}
\includegraphics[width = 0.3\textwidth]{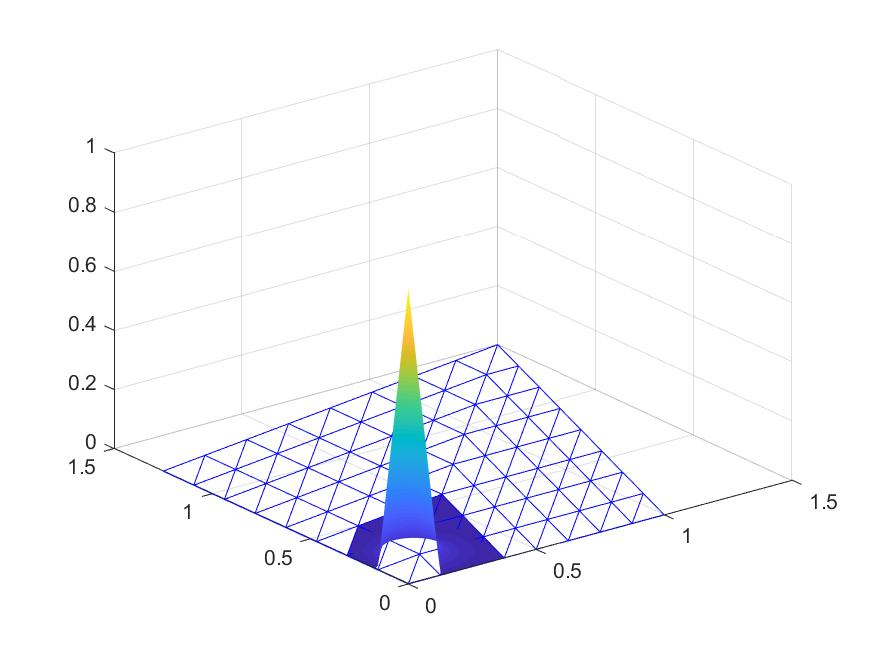}
\includegraphics[width = 0.3\textwidth]{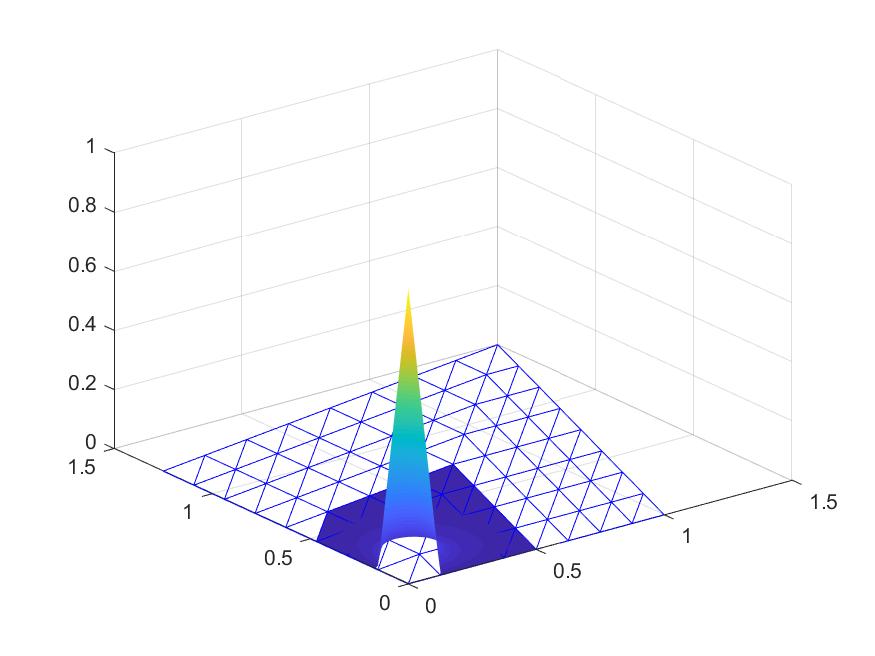}

\includegraphics[width = 0.3\textwidth]{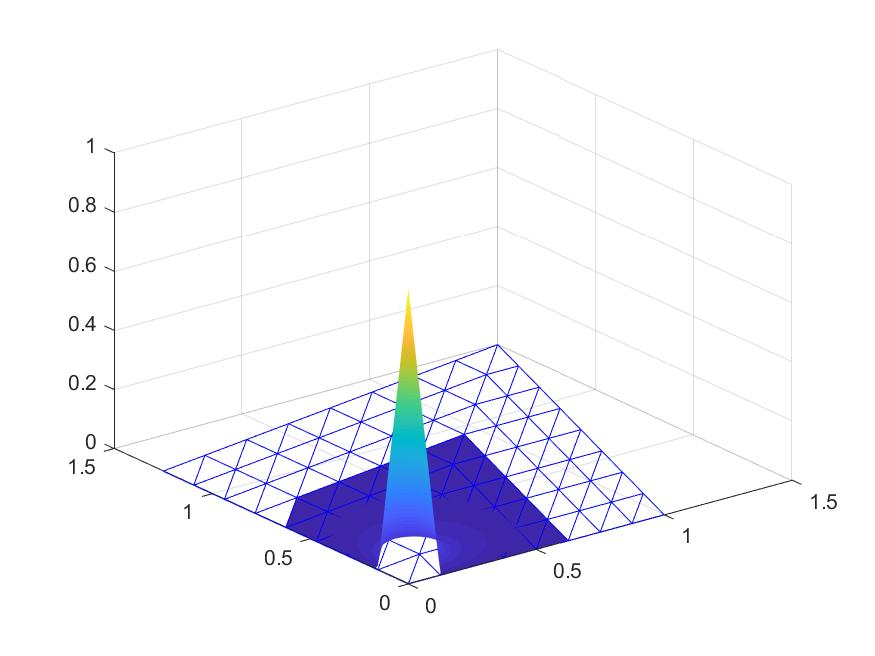}
\includegraphics[width = 0.3\textwidth]{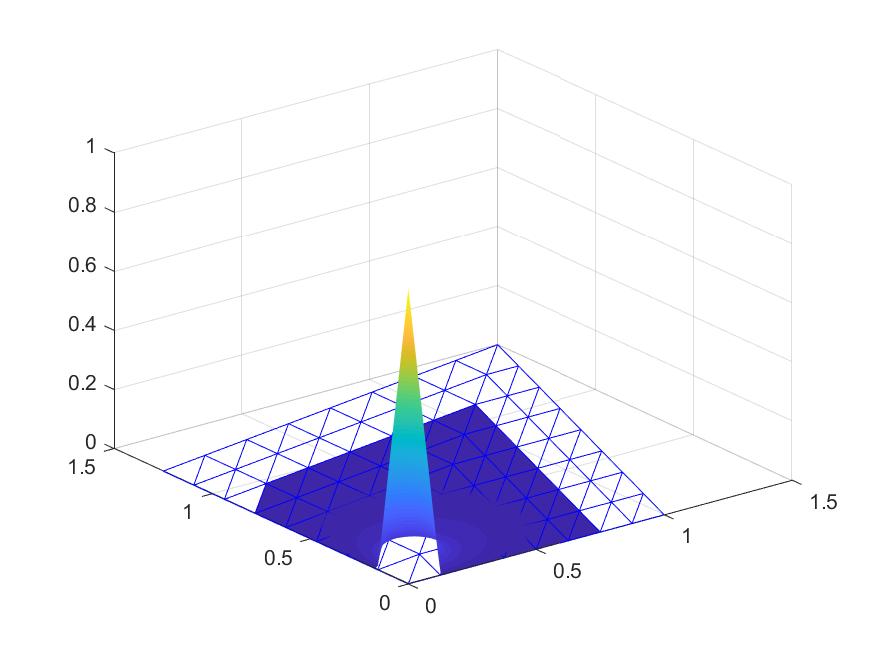}
\includegraphics[width = 0.3\textwidth]{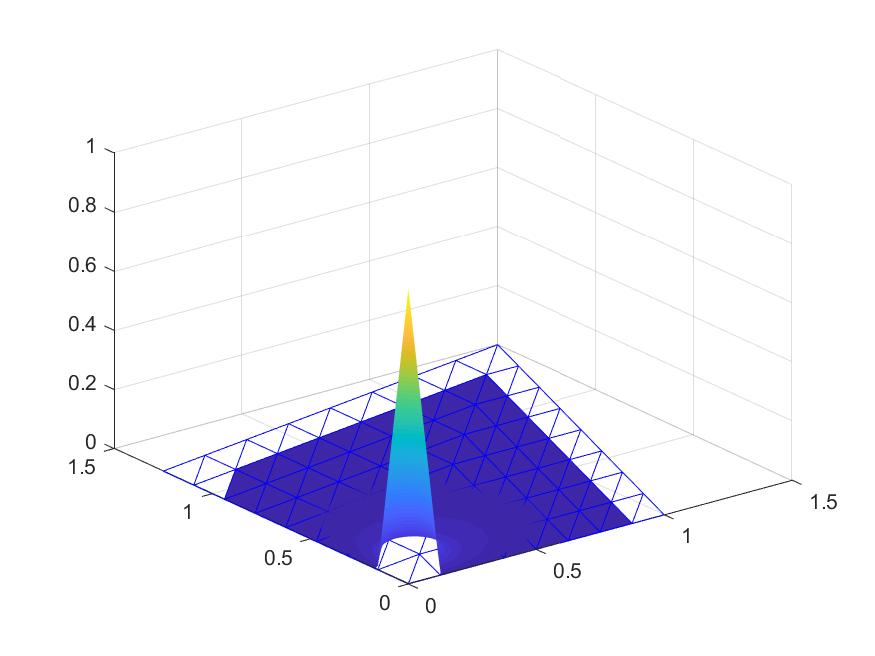}
\caption{A GBC function and its local versions over star$^k(t_0)$ for $k=2, 3, \cdots, 6$. \label{splinegbc2}}
\end{figure}

Let us first show how to compute a boundary GBC locally. Fix $k\ge 1$. 
Let us compute $S_k(\phi_i)$ over the triangulation $\triangle_k=\cup_{t\in T_k} t$ at the boundary vertex $\bfv_i$.  
Then we compare the accuracies of the local GBC functions against the global GBC function supporting at $\bfv_i$.    In
Figure~\ref{splinegbc2}, 
we show the (global) GBC function and its local versions, where 
$t_0:= \bfv_i$.  

Based on those 10,201 equally-spaced points of the bounding box of the polygon 
which fall into the  polygon, we compute 
the maximum errors of these local GBC approximations shown in Figures~\ref{splinegbc2} 
against the global GBC. These maximum errors are presented in Table~\ref{lvg}.
\begin{table}[thpb]
\caption{Local Approximations of $S_i$\label{lvg}}
\centering
\begin{tabular}{|c|c|r||}\hline
no. of rings & max errors & rates \cr \hline
2 & 0.0132 &  \cr
 3 & 0.0072 & 0.5455 \cr
 4 & 0.0042 & 0.5833 \cr
 5 & 0.0024 & 0.5714 \cr
 6 & 0.0011 & 0.4583 \cr\hline
\end{tabular}
\end{table}
From the ratios in Table~\ref{lvg}, we can see that the rate of decay is about $0.5$. 
These numerical results show that we can use local GBCs to replace the (global boundary) 
GBCs for various applications of GBC functions such as 
polygonal domain deformation. This will result a great saving of computational time 
when a polygonal domain is complicated such as the domain of a giraffe or a crocodile.   Indeed, due to the sensitivity of
human eye, an error of $0.0011$ over a graph/image will not be detected by a normal humanbeing.  

Similarly, we can compute local GBC approximations for each interior GBC function 
$R_i= S_{\psi_i}$. For each interior vertex $\bfv_i$, we let star$^k(\bfv_i)$ be the
k-th disk of triangles sharing the vertex $\bfv_i$. 
\begin{figure}[htbp]
\centering
\includegraphics[width = 0.3\textwidth]{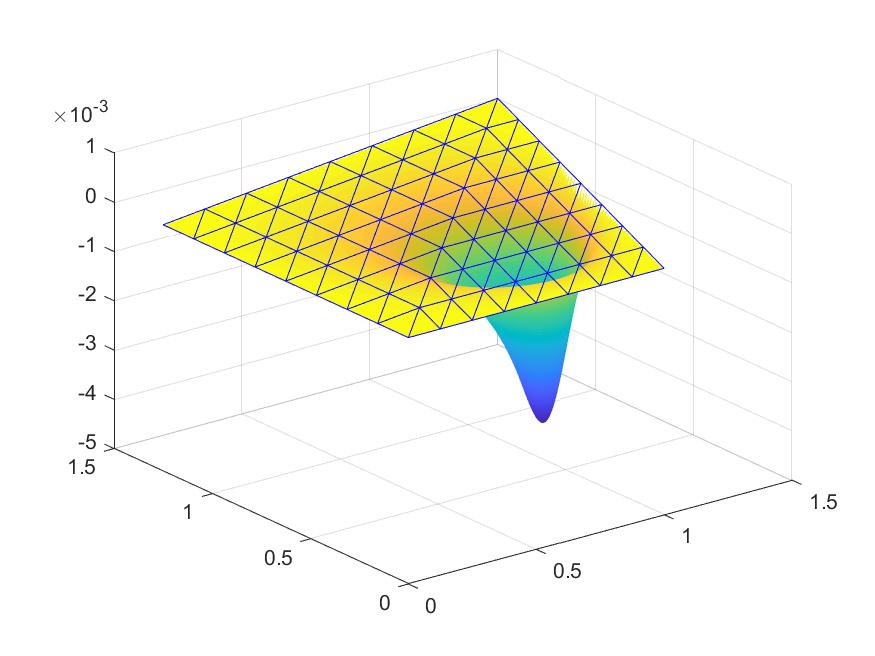}
\includegraphics[width = 0.3\textwidth]{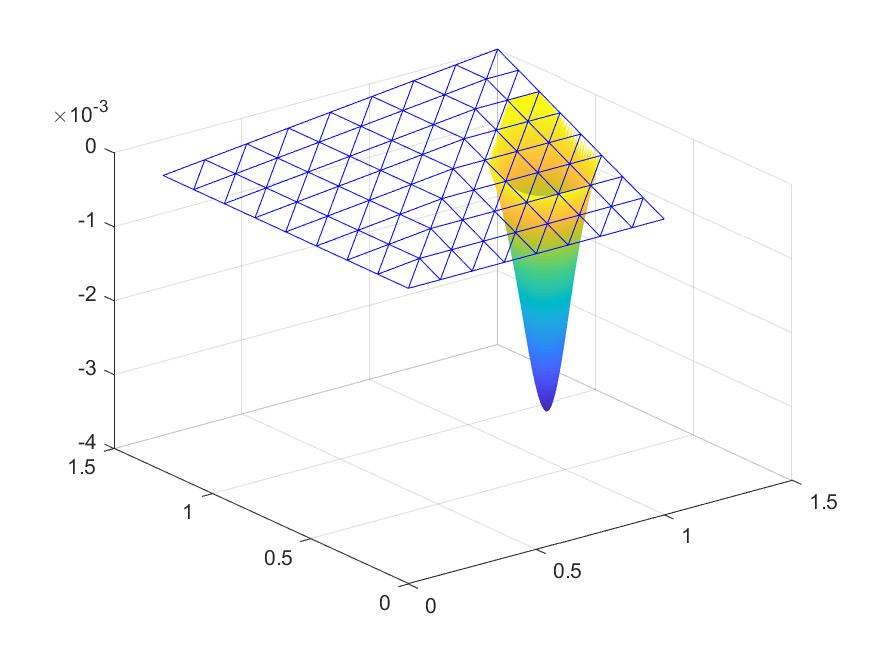}
\includegraphics[width = 0.3\textwidth]{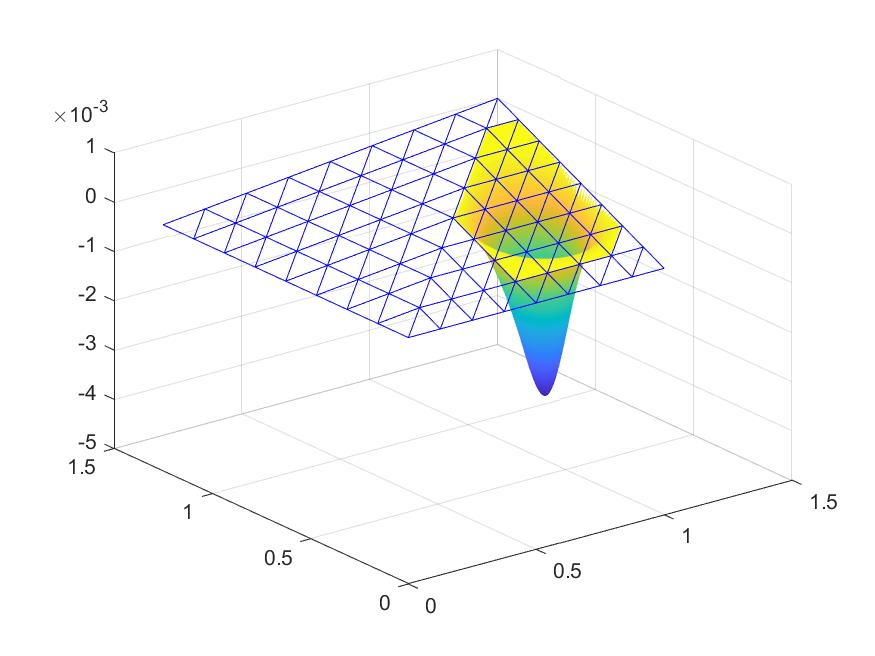}

\includegraphics[width = 0.3\textwidth]{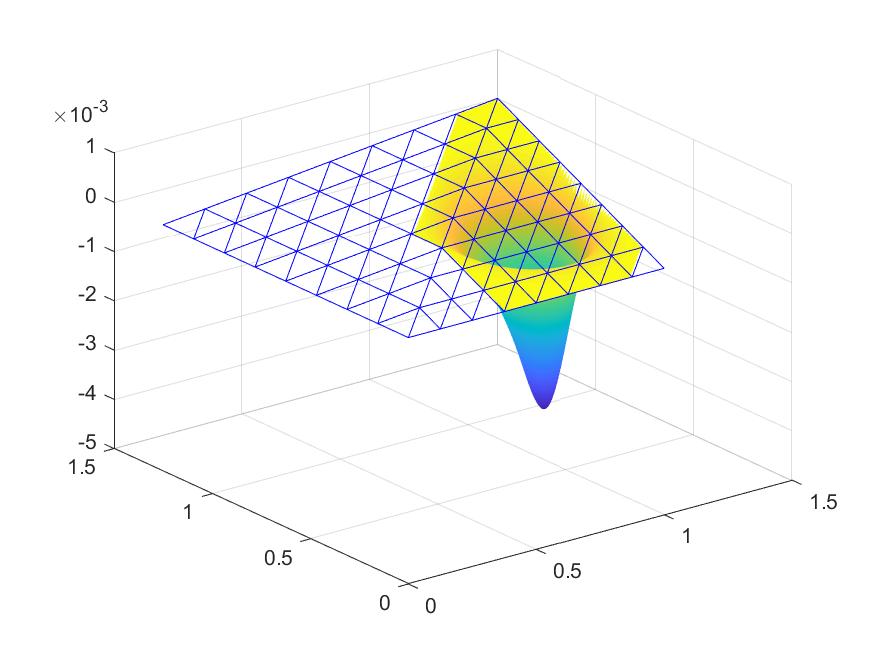}
\includegraphics[width = 0.3\textwidth]{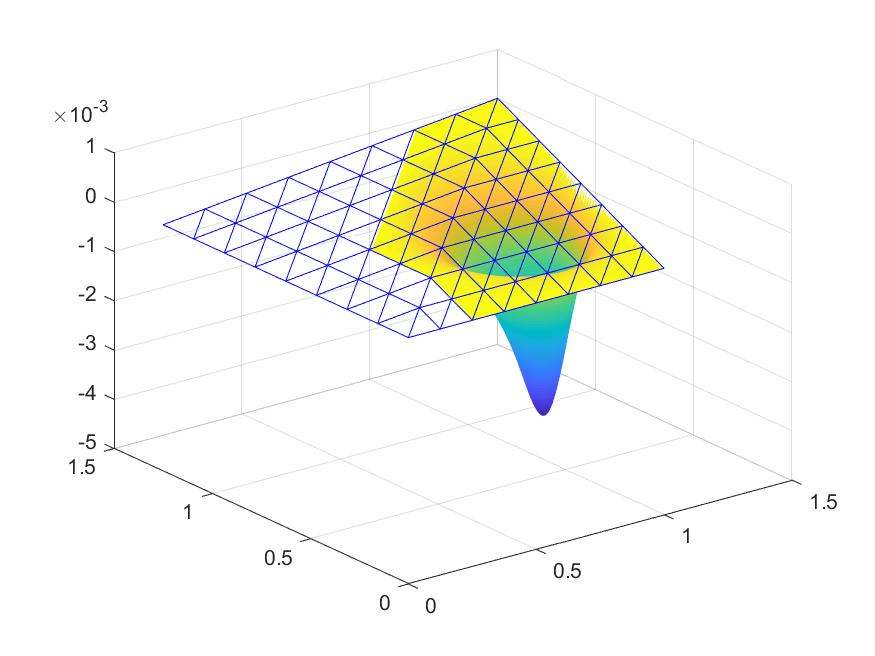}
\includegraphics[width = 0.3\textwidth]{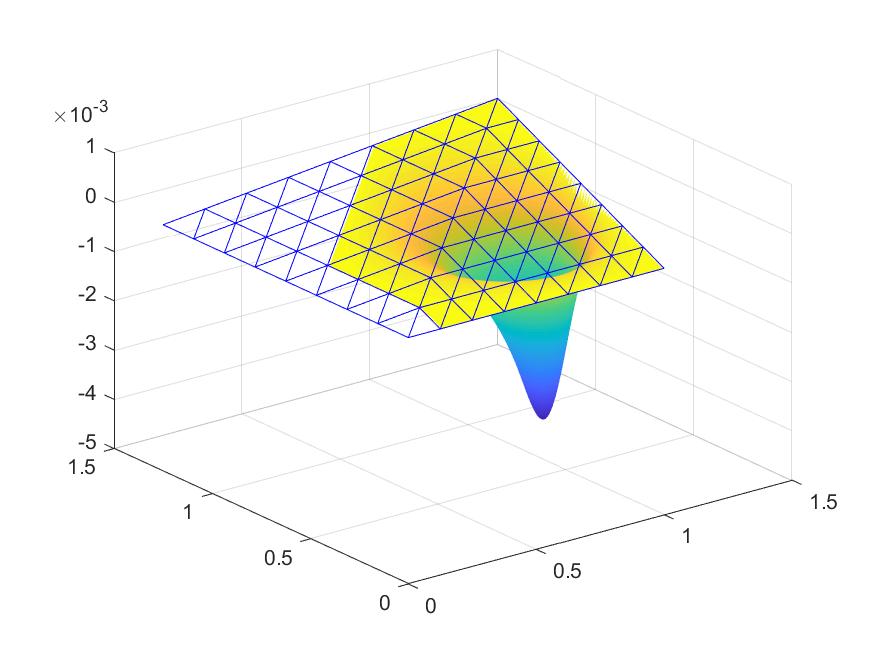}
\caption{An interior GBC function and its local versions 
over star$^k(\bfv_i)$ for $k=4$, $k=5$, and $k=6$. \label{splinegbc4}}
\end{figure}  

The maximum errors of these local GBC approximations shown in Figures~\ref{splinegbc2} 
against the global GBC are given in Table~\ref{lvg2}. Again, we compute the errors based on those 
10,201 equally-spaced points of the bounding box of the polygon which fall into the polygon.  
From the Table~\ref{lvg2}, we can see that the averaged decay rate is about 
$0.68$. 
\begin{table}[thpb]
\caption{Local Approximations of $R_i$\label{lvg2}}
\centering
\begin{tabular}{|c|c|r||}\hline
no. of rings & max errors & rates \cr \hline
 2 & 0.002500 &  \cr
 3 & 0.001800 & 0.72 \cr
 4 & 0.001200 & 0.66 \cr
 5 & 0.000843 & 0.70 \cr
 6 & 0.000559 & 0.66 \cr\hline
\end{tabular}
\end{table}  

Due to the e-locality of these harmonic GBC functions, we can use $k$-local GBC functions to 
approximate the globally GBC functions for an integer $k=2, 3, 4, \ldots$.  

\section{Numerical Approximation of the Poisson Equation}
We now use these boundary-GBC and interior GBC functions to approximate the solution to the 
Poisson equation. In this section, we use a simple domain to demonstrate that these GBC 
functions can indeed be used for numerical solution of PDE.
Besides the domain in the previous section, we also use the convex domain as shown in 
Fig~\ref{cquad}. 

We shall test our method for Poisson equations with various right-hand sides:  
\begin{equation}
\label{Poisson}
\begin{cases}
-\Delta u &=f, \quad \bfx\in \Omega\cr
u &= g, \quad \bfx\in \partial \Omega,
\end{cases}
\end{equation}
where $f$ and $g$ are computed based on testing functions 
\begin{itemize}
\item case 1 $u= 1/(1+x^2+y^2)$, 
\item case 2 $u=x^2+3y^3+4xy$, 
\item case 3 $u=x^4+y^4$, 
\item case 4 $u=\sin(x)\exp(y)$ and 
\item case 5 $u=10\exp(-x^2-y^2)$.
\end{itemize}
Our method is to use a continuous piecewise linear spline based on the
triangulation shown in Figures~\ref{quad} and \ref{cquad} to approximate $f$ and piecewise linear spline on $\partial
\Omega$ to approximate $g$.  Then the GBC approximation of the solution $u$ is given by
\begin{equation}
\label{ourmethod}
u_{gbc} = \sum_{i\in V_B} g(\bfv_i) \phi_i + \sum_{j\in V_I} f(\bfv_j) \psi_j,
\end{equation} 
where $V_B$ and $V_I$ are the index set of the boundary vertices and the index set of interior 
vertices of the triangulation shown in Figure~\ref{quad} or Figure~\ref{cquad}, respectively.  We use bivariate spline 
functions in $S^1_5(\triangle)$ to approximate $\phi_i$ and $\psi_j$. To have a good approximation, we refine the
triangulation shown in Figure~\ref{quad} and Figure~\ref{cquad} once and then compute spline approximation of these
two types of GBC functions. That is, we have
\begin{equation}
\label{ourmethod2}
L_u = \sum_{i\in V_B} g(\bfv_i) S_{\phi_i} + \sum_{j\in V_I} f(\bfv_j) S_{\psi_j},
\end{equation} 

We now show that $L_u$ is a good approximation of $u$. That is, we are now ready to present a 
proof of Theorem~\ref{thm1}. 

\begin{proof}[Proof of Theorem~\ref{thm1}]
First of all, let $s_u$ be the finite element approximation of $u$ over the simplicial partition $\triangle$. 
It is well-known 
that $\|\nabla (u- s_u)\| = O(|\triangle|)$ (cf. \cite{C78}). 
Note that $s_u|_{\partial \Omega}= L_u|_{\partial\Omega}$.  It follows
that 
$$
\|\nabla (u- L_u)\|^2 = \langle \nabla(u- L_u), \nabla (u-s_u)\rangle + \langle \nabla (u- L_u), \nabla( s_u - L_u)\rangle.
$$
We first use Cauchy-Schwarz inequality  to the first term on the right to have
$$
|\langle \nabla(u- L_u), \nabla (u-s_u)\rangle|\le \|\nabla (u- L_u)\| \|\nabla (u-s_u)\| = O(|\triangle|) 
\|\nabla (u- L_u)\|.
$$
Then we use the Green identity to the second term on the right to have
$$
\langle \nabla (u- L_u), \nabla( s_u - L_u)\rangle = -\langle \Delta (u- L_u), s_u- L_u\rangle 
= \langle f- H_f, s_u- u\rangle + \langle f- H_f, u- L_u\rangle
$$
since $\Delta L_g=0$. Note that $H_f$ is a continuous piecewise linear approximation of $f$. Since $H_f$ is a good
approximation of $f$, $\|f- H_f\| = O(|\triangle|^2)$ by the given assumption, we have
$$
|\langle f- H_f, s_u- u\rangle|\le \|f- H_f\| \|s_u- u\|= O(|\triangle|^3)
$$
and
$$
|\langle f- H_f, u- L_u\rangle|\le \|f- H_f\| (\|u\|+\|L_u\|)\le C(|\triangle|^2)
$$
as $\|L_u\|\le \|L_f\|+ \|L_g\|= \|f\|_{L^2(\Omega)} + \|g\|_{L^\infty(\partial\Omega)}+O(|\triangle|^2)$, 
where we have used the GBC property of $S_{\phi_i}$'s and the constant $C$ is dependent on $\|u\|_{L^2(\Omega)}$, 
$\|f\|_{L^2(\Omega)}$, and  $\|g\|_{L^\infty(\partial\Omega)}$.   

Combining the above estimates, we have
$$
\|\nabla (u- L_u)\|^2 \le O(|\triangle|)  \|\nabla (u- L_u)\| +  O(|\triangle|^3) +  C(|\triangle|^2).
$$
In other words, when $|\triangle|$ small enough, we have
$$
(\|\nabla (u- L_u)\|- O(|\triangle|))^2 \le O(|\triangle|^2)+ O(|\triangle|^3) +  C(|\triangle|^2)
$$
or $\|\nabla (u- L_u)\|\le O(|\triangle|)$.  These complete the proof. 
\end{proof}  

\begin{remark}
When the domain $\Omega$ is of uniformly positive reach (cf. \cite{GL20} and \cite{LL22}), we know the solution 
$u$ of the Poisson equation is in $H^2(\Omega)$ when $u|_{\partial\Omega}=0$ and hence, 
\begin{equation}
\label{positivereach}
\|\Delta (u- L_u)\|_{L^2(\Omega)} = \| f - \Delta L_f\|_{L^2(\Omega)}= \| f -\sum_{\bfv_i\in \triangle^{\circ}}f(\bfv_i)
h_i\|_{L^2(\Omega)} = O(|\triangle|).
\end{equation}
Following the argument in \cite{LL22}, there is a positive constant $C>0$ such that 
\begin{equation}
\|u- L_u\|_{H^2(\Omega)} \le C \|\Delta (u- L_u)\|_{L^2(\Omega)},
\end{equation}
where $\|\cdot\|_{H^2(\Omega)}$ is the standard $H^2$ norm in $H^2(\Omega)$.  
\end{remark}

Next we report our numerical results on the approximation of $L_u$ to $u$. 
We use cases 1 through 5 to denote these testing functions mentioned above.  
To see how well our method can approximate the exact solution, we present the maximum errors 
of the GBC solution against the exact solution. In addition,we compare our numerical results 
with the numerical results from the standard finite element method (i.e. continuous linear finite element method).
For simplicity, we compare the maximum error of the numerical solutions from both method against the exact solution, where
the maximum errors are computed based on $1,000\times 1,000$ equally spaced points of the bounding box of the domain 
$\Omega$ which are located inside and on the $\Omega$. See Tables~\ref{result1} and \ref{result2}.   

\begin{table}[thpb]
\caption{Numerical Approximation of $L_u$  over a Nonconvex Quadrilateral \label{result1}}
\centering
\begin{tabular}{|l|c|c|}\hline
& GBCmethod & FEM \cr \hline
case 1  & 0.0179 &  0.0171 \cr
case 2  & 0.3644 &  0.3827 \cr
case 3  & 0.3812 &  0.3506 \cr
case 4  & 0.0610 &   0.0610 \cr
case 5  & 0.1898 &  0.1824 \cr\hline
\end{tabular}
\end{table}

\begin{table}[thpb]
\caption{Numerical Approximation of $L_u$  over a  Convex Quadrilateral \label{result2}}
\centering
\begin{tabular}{|l|c|c|}\hline
& GBCmethod & FEM \cr \hline
case 1 &  0.0188 & 0.0158 \cr
case 2 &  0.2444  & 0.5015 \cr
case 3 &  0.1911  & 0.2063 \cr
case 4 &  0.0175  & 0.0285 \cr
case 5 &  0.2076  & 0.1726 \cr\hline
\end{tabular}
\end{table}

In addition, we refine the triangulation used to generate the numerical approximation 
in Table~\ref{result2} twice 
and compute the maximum errors again. The results are given in Table~\ref{result3}. Note that two methods yield similar
numerical results. For cases 2 and 4, our GBC method is more accurate.  

\begin{table}[thpb]
\caption{Numerical Approximation of $L_u$  over a Triangulation (the first row of each case) and its Refinement (the second
row of each case) \label{result3}}
\centering
\begin{tabular}{|l|c|c|}\hline
& GBCmethod & FEM \cr \hline 
case 1  & 0.00566 &  0.00439 \cr
                  & 0.00155      & 0.00111 \cr \hline 
case 2  & 0.06955 &  0.14491 \cr
                  & 0.01956      & 0.03888  \cr \hline 
case 3  & 0.05576 &  0.05542 \cr 
                  & 0.01524      & 0.01433  \cr \hline 
case 4  & 0.00451 &  0.00737 \cr
                  & 0.00112      & 0.00185  \cr \hline
case 5  & 0.05966 & 0.04446 \cr
                  & 0.01596      &0.01109 \cr \hline
\end{tabular}
\end{table}

From these tables, we can see that the maximum errors of GBC approximation 
$L_u$ are similar to the finite 
element solutions.  These verify Theorem~\ref{thm1} numerically. 
Thus, $L_u$ can indeed be used to approximate the solution of the Poisson equation.  
If the coefficients of these $S_{\phi_i}$ and $S_{\psi_j}$  
can be precomputed and stored, then the 
computational time for numerical solution of the Poisson equation for any right-hand side 
and boundary condition will be greatly reduced. This gives a flexibility of 
solving the Poisson equation, e.g. it allows a modification of a few places in the 
boundary condition and/or a few 
places in the right-hand side to obtain updates of the solution straightforwardly  instead of 
repeatedly solving the system of linear equations again.   
 Furthermore, if each GBC function $S_{\phi_i}$ or $S_{\psi_j}$ is approximated by using its 
k-local versions for small $k$, say $k\le 6$ and all of them are computed individually 
using a GPU simultaneously, the 
computational time for numerical solution of the Poisson equation will be even  reduced 
if the computational 
domain for each $k$-local GBC function is smaller than a quarter of the entire domain 
in the 2D setting.

\section{Conclusions}
In this paper, we define a new type of GBC functions which are called interior-GBC functions 
based on harmonic equations. 
Also, we have demonstrated the e-locality of these GBC functions. That is, we showed that 
each GBC function decays to zero away from its supporting vertex as explained in 
Theorem~\ref{main} and 
Theorem~\ref{main3}.  The proofs are based on a nice result  from \cite{LS09} and 
\cite{GS02} on the decay of locally supported spline functions. 
Furthermore, we showed that the solution to
the Poisson equation can be approximated by using the both GBC functions.  
One can even use the k-local version of
these GBC functions to approximate the solution to Dirichlet problem of the Poisson equation. With a help of GPU 
processes, one can solve the Poisson equation more efficiently in the sense that 
we use a computer with a lot of GPUs to solve the Poisson equation for all $S_i, R_i$. 
Once these solutions are done, for any right-hand side $f$ and boundary condition $g$, 
we simply use the linear combination 
$$
\sum_{i=1}^m f({\bf v}_i)R_i + \sum_{i=1}^N g({\bf v}_i) S_i
$$
to have a numerical solution without solving any system of linear equations any more.

\section{Appendix}
In this Appendix, we mainly justify the approximation of $S_i=S_{\phi_i}$ of 
GBC function $\phi_i$ with supporting
vertex $\bfv_i$ for $i\in V_B$. From (\ref{vPDE}), we know that $\phi_i$ 
satisfies the following weak formulation:
\begin{equation}
\label{weak2}
\langle \nabla \phi_{i,0}, \nabla \psi\rangle = -\langle \nabla G_i, \nabla\psi\rangle, 
\quad \forall 
\psi\in  H^1_0(\Omega),
\end{equation}
where $\phi_{i,0}= \phi_i- G_i$. Together with (\ref{weak}), we have
\begin{equation}
\label{orth}
\langle \nabla (\phi_{i,0}- S_{i,0}), \nabla \psi\rangle=0, \forall \psi\in 
S^r_n(\triangle)\cap H^1_0(\Omega).
\end{equation}
Let $B_{\phi_i}\in S^r_n(\triangle)$ be the best approximation of $\phi_i$ in the spline space 
$S^r_n(\triangle)$ 
satisfying the boundary condition $B_{\phi_i}=G_i$ on $\partial \Omega$. 
It follows that 
\begin{eqnarray}
\|\nabla (\phi_i- S_{\phi_i})\|^2 &=& \langle \nabla (\phi_i- S_{\phi_i}), \nabla (\phi_i- S_{
\phi_i})\rangle\cr
&=& \langle \nabla (\phi_{i,0} - S_{i,0}),  \nabla (\phi_i- B_{\phi_i})\rangle\cr
&\le& \|\nabla (\phi_i- S_{\phi_i})\| \| \nabla (\phi_i- B_{\phi_i})\|
\end{eqnarray}
by using (\ref{orth}). That is, 
$$
\|\nabla (\phi_i- S_{\phi_i})\| \le  \| \nabla (\phi_i- B_{\phi_i})\|.
$$
The computation of the quasi-interpolatory formula in \cite{LS07} shows that we have 
$\|\nabla (\phi_i- B_{\phi_i})\|\le C_{\phi_i} |\triangle|^d$.  Therefore, we have (\ref{cea}). 
Hence, we have (\ref{maxnormest}).

\end{document}